\documentclass{article}
\usepackage{amssymb}

\input{tcilatex}
\begin{document}

\begin{center}
\textbf{The regularity properties and blow-up for convolution wave equations
and applications }

\textbf{Veli\ B. Shakhmurov}

Department of Mechanical Engineering, Istanbul Okan University, Akfirat,
Tuzla 34959 Istanbul, E-mail: veli.sahmurov@okan.edu.tr,

Baku Engineering University

E-mail: veli.sahmurov@gmail.com

\textbf{Rishad Shahmurov}

shahmurov@hotmail.com

University of Alabama Tuscaloosa USA, AL 35487

\textbf{Abstract}
\end{center}

In this paper, the Cauchy problem for linear and nonlinear convolution wave
equations are studied.The equation involves convolution terms with a general
kernel functions whose Fourier transform are operator functions defined in a
Banach space $E$ together with some growth conditions. Here, assuming enough
smoothness on the initial data and the operator functions, the local, global
existence, uniqueness and regularity properties of solutions are established
in terms of fractional powers of given sectorial operator functon.
Furthermore, conditions for finite time blow-up are provided.

By choosing the space $E$ and the operators, the regularity properties the
wide class of nonlocal wave equations in the field of physics are obtained.


\begin{center}
\bigskip\ \ \textbf{AMS: 35Q41, 35L90, 47B25, 35L20, 46E40}

\textbf{1}. \textbf{Introduction}
\end{center}

The aim here, is to study the existence, uniqueness, regularity properties
and blow-up on finite point of solutions to the initial value problem (IVP)
for convolution abstract wave equat\i on (WE) 
\begin{equation}
u_{tt}-a\ast \Delta u+A\ast u=\Delta \left[ g\ast f\left( u\right) \right] 
\text{, }t\in \mathbb{R}_{T}^{n}=\mathbb{R}^{n}\times \left( 0,T\right) , 
\tag{1.1}
\end{equation}%
\begin{equation}
u\left( x,0\right) =\varphi \left( x\right) \text{, }u_{t}\left( x,0\right)
=\psi \left( x\right) \text{ for a.e. }x\in \mathbb{R}^{n},  \tag{1.2}
\end{equation}%
where $A=A\left( x\right) $ is a linear and $g=g\left( x\right) $, $f(u)$
are nonlinear operator functions defined in a Banach space $E$; $a$ is a
complex valued functon on $\mathbb{R}^{n}$, $T\in \left( 0,\right. \left.
\infty \right] $, $\Delta $ denotes the Laplace operator in $\mathbb{R}^{n}$%
, $\varphi \left( x\right) $ and $\psi \left( x\right) $ are the given $E-$%
valued initial functions.

\textbf{Remarke 1.1.} Let $u\in Y^{2,s,p}=W^{2,s,p}\left( \mathbb{R}%
_{T}^{n};E\left( A\right) ,E\right) $, then by J. lions-J. Peetre result
(see e.g. $\left[ \text{27,\ \S 1.8.2}\right] $ the trace operator $%
u\rightarrow \frac{\partial ^{i}u}{\partial t^{i}}\left( x,t\right) $ is
bounded from $Y^{2,s,p}$ to$\ C\left( \mathbb{R}^{n};\left(
Y^{s,p},X_{p}\right) _{\theta _{j},p}\right) $, where 
\[
X_{p}=L^{p}\left( \mathbb{R}^{n};E\right) \text{, }Y^{s,p}=W^{s,p}\left( 
\mathbb{R}^{n};E\left( A\right) ,E\right) \text{, }\theta _{j}=\frac{1+jp}{2p%
}\text{, }j=0\text{, }1, 
\]%
Moreover, if $u\left( x,.\right) \in \left( Y^{s,p},X_{p}\right) _{\theta
_{j},p}$, then under some assumptions that will be stated in the Section 3, $%
f\left( u\right) \in E$ for all $x$, $t\in \mathbb{R}_{T}^{n}$ and the map $%
u\rightarrow f\left( u\right) $ is bounded from $\left(
Y^{s,p},X_{p},\right) _{\frac{1}{2p},p}$ into $E$. Hence, the nonlinear
equation $\left( 1.1\right) $ is satisfied in the Banach space $E$. Here, $%
E\left( A\right) $ denotes a domain of $A$ equipped with graphical norm, $%
\left( Y^{s,p},X_{p}\right) _{\theta ,p}$ is a real interpolation space
between $X_{p}$, $Y^{s,p}$ for $\theta \in \left( 0,1\right) $, $p\in \left[
1,\infty \right] $ (see e.g. $\left[ \text{27,\ \S 1.3}\right] $). The
spaces $X_{p}$, $Y^{s,p}$, $Y^{2,s,p}$ will be defined in Section 1.

The predictions of classical (local) elasticity theory become inaccurate
when the characteristic length of an elasticity problem is comparable to the
atomic length scale. To solution this situation, a nonlocal theory of
elasticity was introduced (see $[$1-3$]$ and the references cited therein)
and the main feature of the new theory is the fact that its predictions were
more down to earth than those of the classical theory. For other
generalizations of elasticity we refer the reader to $[$4-6$]$. The global
existence of the Cauchy problem for Boussinesq type nonlocal equations has
been studied by many authors (see $\left[ \text{11, 14, 21}\right] $ ). Note
that, the existence and uniqueness of solutions and regularity properties
for different type wave equations were considered e.g. in $\left[ \text{4-6}%
\right] $, $\left[ \text{8}\right] $, $\left[ \text{10}\right] $, $\left[ 
\text{17,18}\right] $ and $\left[ \text{31, 32}\right] $. Wave type
equations occur in a wide variety of physical systems, such as in the
propagation of longitudinal deformation waves in an elastic rod,
hydro-dynamical process in plasma, in materials science which describe
spinodal decomposition and in the absence of mechanical stresses (see $\left[
\text{19, 20, 29, 33}\right] $). The $L^{p}$ well-posedness of the Cauchy
problem $(1.1)-\left( 1.2\right) $ depends crucially on the presence of a
suitable kernel. Then the question that naturally arises is which of the
possible forms of the operator functions and kernel functions are relevant
for the global well-posedness of the Cauchy problem $\left( 1.1\right)
-\left( 1.2\right) $. In this study, as a partial answer to this question,
we consider the problem $(1.1)-\left( 1.2\right) $ with a general class of
kernel functions with operator coefficients provide local and global
existence and regularity properties of $(1.1)-\left( 1.2\right) $\ in terms
of fractional powers of operator $A$ in frame of $E-$valued $L^{p}$ spaces.
The kernel functions most frequently used in the literature are particular
cases of this general class of kernel functions in the scalar case, i.e.
when $E=\mathbb{C}$ (here, $\mathbb{C}$ denote the set of complex numbers).
In contrast to the above works, we consider the IVP for nonlocal wave
equation with operator coefficients in $E-$valued function spaces. By
choosing the space $E$, operators $A$ and $g$ in $\left( 1.1\right) -\left(
1.2\right) $, we obtain different classes of nonlocal wave equations which
occur in application. Let we put $E=l_{q}$ and choose $A$, $g$\ as infinite
matrices $\left[ a_{mj}\right] $, $\left[ g_{mj}\right] $, respectively for $%
m,j=1,2,...\infty $. Consider IVP for infinity many system of nonlocal WEs $%
\infty $%
\begin{equation}
\partial _{t}^{2}u_{m}-a\ast \Delta u_{m}+\dsum\limits_{j=1}^{\infty
}a_{mj}\ast u_{m}=  \tag{1.3}
\end{equation}%
\[
\dsum\limits_{j=1}^{\infty }\Delta g_{mj}u_{m}\ast f_{m}\left(
u_{1},u_{2},...,u_{m}\right) \text{, }t\in \left[ 0,T\right] \text{, }x\in 
\mathbb{R}^{n}, 
\]%
\[
u_{m}\left( x,0\right) =\varphi _{m}\left( x\right) \text{, }\partial
_{t}u_{m}\left( x,0\right) =\psi _{m}\left( x\right) \text{, }%
m=1,2,...\infty , 
\]%
where $a_{mj}=a_{mj}\left( x\right) $,$\ g_{mj}=\left( x\right) $ are
complex valued functions, $f_{m}$ are nonlinear functions and $%
u_{j}=u_{j}\left( x,t\right) .$

Then from our results we obtain the existence, uniqueness and regularity
properties of the problem $\left( 1.3\right) $ in terms of fractional powers
of matrix operator $A$ in frame of $l_{q}-$valued $L^{p}$ spaces.

Moreover, let we choose $E=L^{p_{1}}\left( 0,1\right) $ and $A$ to be
degenerated differential operator in $L^{p_{1}}\left( 0,1\right) $ defined
by 
\[
D\left( A\right) =\left\{ u\in W_{\gamma }^{\left[ 2\right] ,p_{1}}\left(
0,1\right) \text{,}\right. \left. \alpha _{k}u^{\left[ \nu _{k}\right]
}\left( 0\right) +\beta _{k}u^{\left[ \nu _{k}\right] }\left( 1\right) =0,%
\text{ }k=1,2\right\} ,\text{ } 
\]%
\begin{equation}
\text{ }A\left( x\right) u=b_{1}\left( x,y\right) u^{\left[ 2\right]
}+b_{2}\left( x,y\right) u^{\left[ 1\right] }\text{, }x\in \mathbb{R}^{n}%
\text{, }y\in \left( 0,1\right) \text{, }\nu _{k}\in \left\{ 0,1\right\} , 
\tag{1.4}
\end{equation}%
\ \ \ where $u^{\left[ i\right] }=\left( y^{\gamma }\frac{d}{dy}\right)
^{\gamma }u$ for $0\leq \gamma <\frac{1}{p_{1}}$, $b_{1}=b_{1}\left(
x,y\right) $ is a cont\i nous, $b_{2}=b_{2}\left( x,y\right) $ is a bounded
functon in $y\in $ $\left[ 0,1\right] $ for a.e. $x\in \mathbb{R}^{n}$, $%
\alpha _{k}$, $\beta _{k}$ are complex numbers and $W_{\gamma }^{\left[ 2%
\right] ,p_{1}}\left( 0,1\right) $ is a weighted Sobolev space defined by 
\[
W_{\gamma }^{\left[ 2\right] ,p_{1}}\left( 0,1\right) =\left\{ {}\right.
u:u\in L^{p_{1}}\left( 0,1\right) \text{, }u^{\left[ 2\right] }\in
L^{p_{1}}\left( 0,1\right) ,\text{ } 
\]%
\[
\left\Vert u\right\Vert _{W_{\gamma }^{\left[ 2\right] ,p_{1}}}=\left\Vert
u\right\Vert _{L^{p_{1}}}+\left\Vert u^{\left[ 2\right] }\right\Vert
_{L^{p_{1}}}<\infty . 
\]%
Then from general results we also obtain the existence, uniqueness and
regularity properties for the nonlocal mixed problem for nonlocal degenerate
PDE 
\begin{equation}
u_{tt}-a\ast \Delta u+\left( b_{1}\frac{\partial ^{\left[ 2\right] }u}{%
\partial y^{2}}+b_{2}\frac{\partial ^{\left[ 1\right] }u}{\partial y}\right)
\ast u=\Delta g\ast f\left( u\right) ,\text{ }  \tag{1.5}
\end{equation}%
\[
x\in \mathbb{R}^{n}\text{, }y\in \left( 0,1\right) \text{, }t\in \left(
0,T\right) \text{, }u=u\left( x,y,t\right) , 
\]%
\ \ \ 

\begin{equation}
\alpha _{ki}u^{\left[ \nu _{k}\right] }\left( x,0,t\right) +\beta _{ki}u^{%
\left[ \nu _{k}\right] }\left( x,1,t\right) =0\text{, }k=1,2,  \tag{1.6}
\end{equation}

\begin{equation}
u\left( x,y,0\right) =\varphi \left( x,y\right) \text{, }u_{t}\left(
x,y,0\right) =\psi \left( x,y\right) \text{.}  \tag{1.7}
\end{equation}

Then from our general results we deduscd the existence, uniqueness and
regularity properties of the problem $\left( 1.5\right) -\left( 1.7\right) $
in terms of fractional powers of operator $A$ defined by $\left( 1.4\right) $
in frame of $L^{p_{1}}\left( 0,1\right) $-valued $L^{p}$ spaces.

It should be noted that, the regularity properties of nonlinear wave
equations in terms of interpolation of spaces are very hard to obtain by the
usual classical methods.

The IVP for abstract hyperbolic equations were studied e.g. in $\left[ \text{%
2}\right] $, $\left[ \text{12}\right] $ and $\left[ \text{22, 23}\right] .$

The strategy is to express the equation $\left( 1.1\right) $ as an integral
equation. To treat the nonlinearity as a small perturbation of the linear
part of the equation, the contraction mapping theorem is used. Also, a
priori estimates on $L^{p}$ norms of solutions of the linearized version are
utilized. The key step is the derivation of the uniform estimate for
solutions of the linearized convolution wave equation. The methods of
harmonic analysis, operator theory, interpolation of Banach spaces and
embedding theorems in Sobolev spaces are the main tools implemented to carry
out the analysis.

In order to state our results precisely, we introduce some notations and
some function spaces.

\begin{center}
\textbf{1.} \textbf{Definitions and} \textbf{Background}
\end{center}

Let $E$ be a Banach space. $L^{p}\left( \Omega ;E\right) $ denotes the space
of strongly measurable $E$-valued functions that are defined on the
measurable subset $\Omega \subset \mathbb{R}^{n}$ with the norm

\[
\left\Vert f\right\Vert _{p}=\left\Vert f\right\Vert _{L^{p}\left( \Omega
;E\right) }=\left( \int\limits_{\Omega }\left\Vert f\left( x\right)
\right\Vert _{E}^{p}dx\right) ^{\frac{1}{p}}\text{, }1\leq p<\infty ,\text{ }
\]

\[
\left\Vert f\right\Vert _{L^{\infty }\left( \Omega ;E\right) }\
=ess\sup\limits_{x\in \Omega }\left\Vert f\left( x\right) \right\Vert _{E}. 
\]

Let $E_{1}$ and $E_{2}$ be two Banach spaces. $\left( E_{1},E_{2}\right)
_{\theta ,p}$ for $\theta \in \left( 0,1\right) $, $p\in \left[ 1,\infty %
\right] $ denotes the real interpolation spaces defined by $K$-method $\left[
\text{27, \S 1.3.2}\right] $. Let $E_{1}$ and $E_{2}$ be two Banach spaces. $%
B\left( E_{1},E_{2}\right) $ will denote the space of all bounded linear
operators from $E_{1}$ to $E_{2}$. For $E_{1}=E_{2}=E$ it will be denoted by 
$B\left( E\right) .$

$\mathbb{N-}$denote the set of natural numbers and $\mathbb{C}$ denotes the
set of complex numbers. Here, \ 
\[
S_{\psi }=\left\{ \lambda \in \mathbb{C}\text{, }\left\vert \arg \lambda
\right\vert \leq \phi \text{, }0\leq \phi <\pi \right\} . 
\]

A closed linear operator\ $A$ is said to be sectorial in a Banach\ space $E$
with bound $M>0$ if $D\left( A\right) $ and $R\left( A\right) $ are dense on 
$E$, $N\left( A\right) =\left\{ 0\right\} $ and 
\[
\left\Vert \left( A+\lambda I\right) ^{-1}\right\Vert _{B\left( E\right)
}\leq M\left\vert \lambda \right\vert ^{-1} 
\]%
for any $\lambda \in S_{\phi }$, $0\leq \phi <\pi $, where $I$ is the
identity operator in $E$, $B\left( E\right) $ is the space of bounded linear
operators in $E;$ $D\left( A\right) $ and $R\left( A\right) $ denote domain
and range of the operator $A.$ It is known that (see e.g.$\left[ \text{27, 
\S 1.15.1}\right] $) there exist the fractional powers\ $A^{\theta }$ of a
sectorial operator $A$. Let $E\left( A^{\theta }\right) $ denote the space $%
D\left( A^{\theta }\right) $ with the graphical norm 
\[
\left\Vert u\right\Vert _{E\left( A^{\theta }\right) }=\left( \left\Vert
u\right\Vert ^{p}+\left\Vert A^{\theta }u\right\Vert ^{p}\right) ^{\frac{1}{p%
}}\text{, }1\leq p<\infty \text{, }0<\theta <\infty . 
\]%
A sectorial operator $A\left( \xi \right) $ for $\xi \in \mathbb{R}^{n}$ is
said to be uniformly sectorial in a Banach space $E$, if $D\left( A\left(
\xi \right) \right) $ is independent of $\xi $ and the uniform estimate 
\[
\left\Vert \left( A+\lambda I\right) ^{-1}\right\Vert _{B\left( E\right)
}\leq M\left\vert \lambda \right\vert ^{-1} 
\]%
holds for any $\lambda \in S_{\phi }.$

A linear operator\ $A=A\left( \xi \right) $ belongs to $\sigma \left(
M_{0},\omega ,E\right) $ (see $\left[ \text{23}\right] $ \S\ 11.2) if $%
D\left( A\right) $ is dense on $E,$ $D\left( A\left( \xi \right) \right) $
is independent of $\xi \in \mathbb{R}^{n}$ and for $\func{Re}\lambda >\omega 
$ the uniform estimate holds 
\[
\left\Vert \left( A\left( \xi \right) -\lambda ^{2}I\right) ^{-1}\right\Vert
_{B\left( E\right) }\leq M_{0}\left\vert \func{Re}\lambda -\omega
\right\vert ^{-1}\text{. } 
\]

\textbf{Remark 1.1. }It is known (see e.g. $\left[ \text{22, \S\ 1.6}\right] 
$, Theorem 6.3)\ that if $A\in \sigma \left( M_{0},\omega ,E\right) $ and $%
0\leq \alpha <1$ then it is generates a bounded group operator $U_{A}\left(
t\right) $ satisfying 
\begin{equation}
\left\Vert U_{A}\left( t\right) \right\Vert _{B\left( E\right) }\leq
Me^{\omega \left\vert t\right\vert }\text{, }\left\Vert A^{\alpha
}U_{A}\left( t\right) \right\Vert _{B\left( E\right) }\leq M\left\vert
t\right\vert ^{-\alpha }\text{, }t\in \left[ 0,T\right] .  \tag{2.1}
\end{equation}

Let $E$ be a Banach space. $S=S(\mathbb{R}^{n};E)$ denotes $E$-valued
Schwartz class, i.e. the space of all $E-$valued rapidly decreasing smooth
functions on $\mathbb{R}^{n}$ equipped with its usual topology generated by
seminorms. $S(\mathbb{R}^{n};\mathbb{C})$ denoted by $S$. Let $S^{\prime }(%
\mathbb{R}^{n};E)$ denote the space of all continuous linear functions from $%
S$ into $E$, equipped with the bounded convergence topology. Recall $S(%
\mathbb{R}^{n};E)$ is norm dense in $L^{p}(\mathbb{R}^{n};E)$ when $1\leq
p<\infty $. Let $m$ be a positive integer. $W^{m,p}\left( \Omega ;E\right) $
denotes an $E-$valued Sobolev space of all functions $u\in L^{p}\left(
\Omega ;E\right) $ that have the generalized derivatives $\frac{\partial
^{m}u}{\partial x_{k}^{m}}\in L^{p}\left( \Omega ;E\right) $ with the norm 
\[
\ \left\Vert u\right\Vert _{W^{m,p}\left( \Omega ;E\right) }=\left\Vert
u\right\Vert _{L^{p}\left( \Omega ;E\right)
}+\sum\limits_{k=1}^{n}\left\Vert \frac{\partial ^{m}u}{\partial x_{k}^{m}}%
\right\Vert _{L^{p}\left( \Omega ;E\right) }<\infty . 
\]%
\ \ $\ \ $

Let $W^{s,p}\left( \mathbb{R}^{n};E\right) $ denotes the fractional Sobolev
space of order $s\in \mathbb{R}$, that is defined as: 
\[
W^{s,p}\left( E\right) =W^{s,p}\left( \mathbb{R}^{n};E\right) =\left\{ u\in
\right. S^{\prime }(\mathbb{R}^{n};E), 
\]%
\[
\left. \left\Vert u\right\Vert _{W^{s,p}\left( E\right) }=\left\Vert \mathbb{%
F}^{-1}\left( I+\left\vert \xi \right\vert ^{2}\right) ^{\frac{s}{2}}\hat{u}%
\right\Vert _{L^{p}\left( \mathbb{R}^{n};E\right) }<\infty \right\} . 
\]%
It clear that $W^{0,p}\left( \mathbb{R}^{n};E\right) =L^{p}\left( \mathbb{R}%
^{n};E\right) $. Let $E_{0}$ and $E$ be two Banach spaces and $E_{0}$ is
continuously and densely embedded into $E$. Here, $W^{s,p}\left( \mathbb{R}%
^{n};E_{0},E\right) $ denote the Sobolev-Lions type space i.e., 
\[
W^{s,p}\left( \mathbb{R}^{n};E_{0},E\right) =\left\{ u\in W^{s,p}\left( 
\mathbb{R}^{n};E\right) \cap L^{p}\left( \mathbb{R}^{n};E_{0}\right)
,\right. \text{ } 
\]%
\[
\left. \left\Vert u\right\Vert _{W^{s,p}\left( \mathbb{R}^{n};E_{0},E\right)
}=\left\Vert u\right\Vert _{L^{p}\left( \mathbb{R}^{n};E_{0}\right)
}+\left\Vert u\right\Vert _{W^{s,p}\left( \mathbb{R}^{n};E\right) }<\infty
\right\} . 
\]

In a similar way, we define the following Sobolev-Lions type space:

\[
W^{2,s,p}\left( \mathbb{R}_{T}^{n};E_{0},E\right) =\left\{ {}\right. u\in
S^{\prime }\left( \mathbb{R}_{T}^{n};E_{0}\right) \text{, }\partial
_{t}^{2}u\in L^{p}\left( \mathbb{R}_{T}^{n};E\right) , 
\]%
\[
\mathbb{F}_{x}^{-1}\left( I+\left\vert \xi \right\vert ^{2}\right) ^{\frac{s%
}{2}}\hat{u}\in L^{p}\left( \mathbb{R}_{T}^{n};E\right) \text{, }\left\Vert
u\right\Vert _{W^{2,s,p}\left( \mathbb{R}_{T}^{n};E_{0},E\right) }= 
\]%
\[
\left\Vert \partial _{t}^{2}u\right\Vert _{L^{p}\left( \mathbb{R}%
_{T}^{n};E\right) }+\left. \left\Vert \mathbb{F}_{x}^{-1}\left( I+\left\vert
\xi \right\vert ^{2}\right) ^{\frac{s}{2}}\hat{u}\right\Vert _{L^{p}\left( 
\mathbb{R}_{T}^{n};E\right) }<\infty \right\} . 
\]

Let $1\leq p\leq q<\infty $. A function $\Psi \in L^{\infty }(\mathbb{R}%
^{n}) $ is called a Fourier multiplier from $L^{p}(\mathbb{R}^{n};E)$ to $%
L^{q}(\mathbb{R}^{n};E)$ if the map $P$: $u\rightarrow \mathbb{F}^{-1}\Psi
(\xi )\mathbb{F}u$ for $u\in S(\mathbb{R}^{n};E)$ is well defined and
extends to a bounded linear operator

\[
P\text{: }L^{p}(\mathbb{R}^{n};E)\rightarrow L^{q}(\mathbb{R}^{n};E). 
\]

A Banach space $E$ is called a UMD space\ if the Hilbert operator

\[
(Hf)(x)=\lim\limits_{\varepsilon \rightarrow 0}\int\limits_{\left\vert
x-y\right\vert >\varepsilon }\frac{f(y)}{x-y}dy 
\]%
is initially defined on $S(\mathbb{R};E)$ and is bounded in $L^{p}(\mathbb{R}%
;E)$, $p\in (1,\infty )$. UMD spaces include e.g. $L_{p},l_{p}$ spaces and
Lorentz spaces $L_{pq}$, $p,q\in (1,\infty )$ (see e.g $\left[ \text{13}%
\right] $).

A set $K\subset B(E_{1},E_{2})$ is called $R-$bounded (see e.g $[$13$]$) if
there is a constant $C>0$ such that for all $T_{1},T_{2},...,T_{m}\in K$ and 
$u_{1},u_{2},...,u_{m}\in E_{1},$ $m\in \mathbb{N}$

\[
\int\limits_{0}^{1}\left\Vert
\sum\limits_{j=1}^{m}r_{j}(y)T_{j}u_{j}\right\Vert _{E_{2}}dy\leq
C\int\limits_{0}^{1}\left\Vert \sum\limits_{j=1}^{m}r_{j}(y)u_{j}\right\Vert
_{E_{1}}dy, 
\]%
where $\left\{ r_{j}\right\} $ is a sequence of independent symmetric $%
\left\{ -1;1\right\} -$valued random variables on $\left[ 0,1\right] $. The
smallest $C$ for which the above estimate holds is called the $R-$bound of $%
K $ and denoted by $R\left( K\right) .$

Note that, in Hilbert spaces every norm bounded set is $R-$bounded.
Therefore, all sectorial operators are $R-$sectorial in Hilbert spaces.

Sometimes we use one and the same symbol $C$ without distinction in order to
denote positive constants which may differ from each other even in a single
context. When we want to specify the dependence of such a constant on a
parameter, say $\alpha $, we write $C_{\alpha }$. Moreover, for $u$, $%
\upsilon >0$ the relations $u\lesssim \upsilon $, $u$\ $\approx $ $\upsilon $
means that there exist positive constants $C,$ $C_{1},$ $C_{2}$ \
independent on $u$ and $\upsilon $ such that, respectively 
\[
u\leq C\upsilon \text{, }C_{1}\upsilon \leq u\leq C_{2}\upsilon . 
\]

The paper is organized as follows: In Section 1, some definitions and
background are given. In Section 2, we obtain the existence of unique
solution and a priory estimates for solution of the linearized problem $%
(1.1)-\left( 1.2\right) $. In Section 3, we show the existence and
uniqueness of local strong solution of the problem $(1.1)-\left( 1.2\right) $%
. In the Section 4, we show the same applications of the problem $%
(1.1)-\left( 1.2\right) .$

Sometimes we use one and the same symbol $C$ without distinction in order to
denote positive constants which may differ from each other even in a single
context. When we want to specify the dependence of such a constant on a
parameter, say $h$, we write $C_{h}$.

\begin{center}
\textbf{2. Estimates for linearized equation}
\end{center}

In this section, we make the necessary estimates for solutions of the Cauchy
problem for the convolution linear WE 
\begin{equation}
u_{tt}-a\ast \Delta u+A\ast u=g\left( x,t\right) \text{, }x\in \mathbb{R}^{n}%
\text{, }t\in \left( 0,T\right) \text{, }T\in \left( 0,\right. \left. \infty %
\right] ,  \tag{2.1}
\end{equation}%
\begin{equation}
u\left( x,0\right) =\varphi \left( x\right) \text{, }u_{t}\left( x,0\right)
=\psi \left( x\right) \text{ for a.e. }x\in \mathbb{R}^{n},  \tag{2.2}
\end{equation}%
where $A=A\left( x\right) $ is a linear operator function in a Banach space $%
E$ and $a\geq 0.$

Let $A$ be a sectorial operator in $E$. Here, 
\[
X_{p}=L^{p}\left( \mathbb{R}^{n};E\right) \text{, }X_{p}\left( A^{\gamma
}\right) =L^{p}\left( \mathbb{R}^{n};E\left( A^{\gamma }\right) \right) 
\text{, }1\leq p,\text{ }q\leq \infty ,
\]%
\[
Y\text{ }^{s,p}=Y\text{ }^{s,p}\left( E\right) =W^{s,p}\left( \mathbb{R}%
^{n};E\right) \text{, }Y_{q}^{s,p}\left( E\right) =Y\text{ }^{s,p}\left(
E\right) \cap X_{q}\text{, }
\]%
\[
\left\Vert u\right\Vert _{Y_{q}^{s,p}}=\left\Vert u\right\Vert
_{W^{s,p}\left( \mathbb{R}^{n};E\right) }+\left\Vert u\right\Vert
_{X_{q}}<\infty \text{,}
\]%
\[
W^{s,p}\left( A^{\gamma }\right) =W^{s,p}\left( \mathbb{R}^{n};E\left(
A^{\gamma }\right) \right) \text{, }0<\gamma \leq 1,
\]%
\[
Y^{s,p}=Y^{s,p}\left( A,E\right) =W^{s,p}\left( \mathbb{R}^{n};E\left(
A\right) ,E\right) \text{, }Y^{2,s,p}=Y^{2,s,p}\left( A,E\right) =
\]%
\[
W^{2,s,p}\left( \mathbb{R}_{T}^{n};E\left( A\right) ,E\right) \text{, }%
Y_{q}^{s,p}\left( A;E\right) =Y^{s,p}\left( E\right) \cap X_{q}\left(
A\right) ,
\]%
\[
\left\Vert u\right\Vert _{Y_{q}^{s,p}\left( A,E\right) }=\left\Vert
u\right\Vert _{Y^{s,p}\left( E\right) }+\left\Vert u\right\Vert
_{X_{q}\left( A\right) }<\infty \text{, }
\]%
\[
\text{ }\mathbb{E}_{0p}=\left( Y^{s,p}\left( A,E\right) ,X_{p}\right) _{%
\frac{1}{2p},p}\text{, }\mathbb{E}_{1p}=\left( Y^{s,p}\left( A,E\right)
,X_{p}\right) _{\frac{1+p}{2p},p}.
\]%
Let $\hat{X}_{p}\left( A^{\alpha }\right) $ denotes the $D\left( A^{\alpha
}\right) $-value function space with norm%
\[
\left\Vert u\right\Vert _{\hat{X}_{p}\left( A^{\alpha }\right) }=\left\Vert
A^{\alpha }\ast u\right\Vert _{X_{p}}.
\]%
Let 
\[
Y_{0}\left( A^{\alpha }\right) =\mathbb{E}_{0p}\cap \hat{X}_{1}\left(
A^{\alpha }\right) \text{, }Y_{1}\left( A^{\alpha }\right) =\mathbb{E}%
_{1p}\cap \hat{X}_{1}\left( A^{\alpha }\right) \text{.}
\]%
\textbf{Remark 2.1.} By properties of real interpolation of Banach spaces
and interpolation of the intersection of the spaces (see e.g. $\left[ \text{%
27, \S 1.3}\right] $) we obtain 
\[
\text{ }\mathbb{E}_{0p}=\left( Y^{s,p}\left( A,E\right) \cap
X_{p},X_{p}\right) _{\frac{1}{2p},p}=\left( Y^{s,p}\left( E\right)
,X_{p}\right) _{\frac{1}{2p},p}\cap \left( X_{p}\left( A\right)
,X_{p}\right) _{\frac{1}{2p},p}=
\]%
\[
W^{s\left( 1-\frac{1}{2p}\right) ,p}\left( \mathbb{R}^{n};E\right) \cap
L^{p}\left( \mathbb{R}^{n};\left( E\left( A\right) ,E\right) _{\frac{1}{2p}%
,p}\right) =
\]%
\[
W^{s\left( 1-\frac{1}{2p}\right) ,p}\left( \mathbb{R}^{n};\left( E\left(
A\right) ,E\right) _{\frac{1}{2p},p},E\right) .
\]%
In a similar way, we have%
\[
\mathbb{E}_{1p}=\left( Y^{s,p}\left( A,E\right) \cap X_{p},X_{p}\right) _{%
\frac{1+p}{2p},p}=W^{\frac{s\left( p-1\right) }{2p},p}\left( \mathbb{R}%
^{n};\left( E\left( A\right) ,E\right) _{\frac{1+p}{2p},p},E\right) .
\]

\textbf{Remark 2.2. }Let $A$ be a sectorial operator in a Banach space $E$.
In view of interpolation by domain of sectorial operators (see e.g.$\left[ 
\text{27,\ \S 1.8.2}\right] $) we have the following relation 
\[
E\left( A^{1-\theta +\varepsilon }\right) \subset \left( E\left( A\right)
,E\right) _{\theta ,p}\subset E\left( A^{1-\theta -\varepsilon }\right) 
\]%
for $0<\theta <$ $1$ and$\ 0<\varepsilon <$ $1-\theta .$

Note that from J. lions-J. Peetre result (see e.g. $\left[ \text{27,\ \S %
1.8.2}\right] $ we obtan the following result.

\textbf{Lemma A}$_{1}$. The trace operator $u\rightarrow \frac{\partial ^{i}u%
}{\partial t^{i}}\left( x,t\right) $ is bounded from $Y^{2,s,p}\left(
A,E\right) $ into$\ $%
\[
C\left( \mathbb{R}^{n};\left( Y^{s,p}\left( A,E\right) ,X_{p}\right)
_{\theta _{j},p}\right) \text{, }\theta _{j}=\frac{1+jp}{2p}\text{, }j=0,1. 
\]

Let $\hat{A}\left( \xi \right) $\ be the Fourier transformation of $A\left(
x\right) $, i.e. $\hat{A}\left( \xi \right) =\mathbb{F}\left( A\left(
x\right) \right) $. We assume that $\hat{A}\left( \xi \right) $ is uniformly
sectorial operator in a Banach\ space $E$. Let%
\[
\eta =\eta \left( \xi \right) =\left[ \hat{a}\left( \xi \right) \left\vert
\xi \right\vert ^{2}+\hat{A}\left( \xi \right) \right] ^{\frac{1}{2}}. 
\]

Let $A$ be a generator of a strongly continuous cosine operator function in
a Banach space $E$ defined by formula%
\[
C\left( t\right) =\frac{1}{2}\left( e^{itA^{\frac{1}{2}}}+e^{-itA^{\frac{1}{2%
}}}\right) 
\]%
(see e.g. $\left[ \text{12, \S 11}\right] $). Then, from the definition of
sine operator-function $S\left( t\right) $ we have%
\[
S\left( t\right) u=\dint\limits_{0}^{t}C\left( \sigma \right) ud\sigma \text{%
, i.e. }S\left( t\right) u=\frac{1}{2i}A^{-\frac{1}{2}}\left( e^{itA^{\frac{1%
}{2}}}-e^{-itA^{\frac{1}{2}}}\right) . 
\]

Let 
\begin{equation}
\eta _{\pm }\left( \xi \right) =e^{it\eta \left( \xi \right) }\pm e^{-it\eta
\left( \xi \right) }\text{, }C\left( t\right) =C\left( \xi ,t\right) =\frac{%
\eta _{+}\left( \xi \right) }{2},  \tag{2.3}
\end{equation}

\[
\text{ }S\left( t\right) =S\left( \xi ,t\right) =\eta ^{-1}\left( \xi
\right) \frac{\eta _{-}\left( \xi \right) }{2i}.
\]%
\textbf{Condition 2.1. }Assume: (1) $\eta \left( \xi \right) \neq 0$\ for
all $\xi \in \mathbb{R}^{n}$; (2) $\hat{a}\in C^{\left( m\right) }\left( 
\mathbb{R}^{n}\right) $ such that 
\begin{equation}
\hat{a}\left( \xi \right) \left\vert \xi \right\vert ^{2}\in S_{\varphi _{1}}%
\text{, }\left( 1+\left\vert \xi \right\vert ^{2}\right) ^{-\left( \frac{s}{2%
}-2\right) }\left\vert D^{\beta }\hat{a}\left( \xi \right) \right\vert \leq
C_{0}\text{, }  \tag{2.0}
\end{equation}%
\[
m=\left\vert \beta \right\vert >1+\frac{n}{p}\text{, }p\in \left( 1,\infty
\right) \text{ for all }\xi \in \mathbb{R}^{n};
\]
(3) $\hat{A}\left( \xi \right) $ is an uniformly $R-$sectorial operator in
UMD space $E$ such that $\hat{A}\left( \xi \right) \in \sigma \left(
M_{0},\omega ,E\right) $; (4)$\ \hat{A}\left( \xi \right) $ is a
differentiable operator function with independent of $\xi $ domain $D\left(
D^{\beta }\hat{A}\left( \xi \right) \right) =D\left( \hat{A}\right) =$ $%
D\left( A\right) $ for $\beta =\left( \beta _{1},\beta _{2},...,\beta
_{n}\right) $ and $\left\vert \beta \right\vert \leq n$. Moreover, the
following uniform estimate holds 
\[
\left\Vert \left[ D^{\beta }\hat{A}\left( \xi \right) \right] \eta ^{-\gamma
}\left( \xi \right) \right\Vert _{B\left( E\right) }\leq M\text{ for }%
0<\gamma <1-\frac{1}{2p};
\]

(5) $\varphi \in $ $\mathbb{E}_{0p}$ and $\psi \in $ $\mathbb{E}_{1p}$.

First we need the following lemmas:

\textbf{Lemma 2.1. }Let the assumption (1) of Condition 2.1 holds. Then,
problem $\left( 2.1\right) -\left( 2.2\right) $ has a unique solution.

\textbf{Proof. }By using of the Fourier transform, we get from $(2.1)-\left(
2.2\right) $:%
\begin{equation}
\hat{u}_{tt}\left( \xi ,t\right) +\eta ^{2}\left( \xi \right) \hat{u}\left(
\xi ,t\right) =\hat{g}\left( \xi ,t\right) ,\text{ }  \tag{2.4}
\end{equation}%
\[
\hat{u}\left( \xi ,0\right) =\hat{\varphi}\left( \xi \right) \text{, }\hat{u}%
_{t}\left( \xi ,0\right) =\hat{\psi}\left( \xi \right) , 
\]%
where $\hat{u}\left( \xi ,t\right) $ is a Fourier transform of $u\left(
x,t\right) $ in $x$ and $\hat{\varphi}\left( \xi \right) $, $\hat{\psi}%
\left( \xi \right) $ are Fourier transform of $\varphi $ and $\psi $,
respectively. By virtue of $\left[ \text{12, \S\ 11.2,4}\right] $ we obtain
that $\eta \left( \xi \right) $ is a generator of a strongly continuous
cosine operator function and problem $(2.4)$ has a unique solution for all $%
\xi \in \mathbb{R}^{n}$ exspressing as%
\begin{equation}
\hat{u}\left( \xi ,t\right) =C\left( \xi ,t\right) \hat{\varphi}\left( \xi
\right) +S\left( \xi ,t\right) \hat{\psi}\left( \xi \right)
+\dint\limits_{0}^{t}S\left( \xi ,t-\tau \right) \hat{g}\left( \xi ,\tau
\right) d\tau ,  \tag{2.5}
\end{equation}%
i.e. problem $(2.1)-\left( 2.2\right) $ has a unique solution 
\begin{equation}
u\left( x,t\right) =C_{1}\left( t\right) \varphi +S_{1}\left( t\right) \psi
+Qg,  \tag{2.6}
\end{equation}%
where $C_{1}\left( t\right) $, $S_{1}\left( t\right) $, $Q$ are linear
operator functions defined by 
\[
C_{1}\left( t\right) \varphi =\mathbb{F}^{-1}\left[ C\left( \xi ,t\right) 
\hat{\varphi}\left( \xi \right) \right] \text{, }S_{1}\left( t\right) \psi =%
\mathbb{F}^{-1}\left[ S\left( \xi ,t\right) \hat{\psi}\left( \xi \right) %
\right] , 
\]

\[
Qg=\mathbb{F}^{-1}\tilde{Q}\left( \xi ,t\right) \text{, }\tilde{Q}\left( \xi
,t\right) =\dint\limits_{0}^{t}\mathbb{F}^{-1}\left[ S\left( \xi ,t-\tau
\right) \hat{g}\left( \xi ,\tau \right) \right] d\tau . 
\]%
\textbf{Theorem 2.1. }Assume the Condition 2.1 holds and $s>1+\frac{n}{p}$
with $p\in \left( 1,\infty \right) $. Let $0<\alpha <1-\frac{1}{2p}$. Then
for $\varphi \in Y_{0}\left( A^{\alpha }\right) $, $\psi \in $ $Y_{1}\left(
A^{\alpha }\right) $ and $g\left( x,t\right) \in Y_{1}^{s,p}$ problem $%
(2.1)-(2.2)$ has a unique generalized solution $u(x,t)\in C^{2}\left( \left[
0,T\right] ;X_{\infty }\right) $. Moreover, the following estimate holds 
\begin{equation}
\left\Vert A^{\alpha }\ast u\right\Vert _{X_{\infty }}+\left\Vert A^{\alpha
}\ast u_{t}\right\Vert _{X_{\infty }}\leq C_{0}\left[ \left\Vert \varphi
\right\Vert _{Y_{01}^{\alpha }\left( A\right) }\right. +  \tag{2.7}
\end{equation}

\[
\left\Vert \psi \right\Vert _{Y_{11}^{\alpha }\left( A\right) }+\left.
\dint\limits_{0}^{t}\left( \left\Vert g\left( .,\tau \right) \right\Vert
_{Y_{1}^{s,p}}+\left\Vert g\left( .,\tau \right) \right\Vert _{X_{1}}\right)
d\tau \right] , 
\]%
uniformly in $t\in \left[ 0,T\right] $, where the constant $C_{0}>0$ depends
only on $A$, the space $E$ and initial data.

\textbf{Proof. }From Lemma 2.1 we obtain that problem $(2.1)-(2.2)$ has a
unique generalized solution $u(x,t)\in C^{2}\left( \left[ 0,T\right]
;Y^{s,p}\left( A;E\right) \right) $ for $\varphi \in \mathbb{E}_{0p}$, $\psi
\in \mathbb{E}_{1p}$ and $g\left( .,t\right) \in Y_{1}^{s,p}$. Let $N\in 
\mathbb{N}$ and 
\[
\Pi _{N}=\left\{ \xi :\xi \in \mathbb{R}^{n}\text{, }\left\vert \xi
\right\vert \leq N\right\} \text{, }\Pi _{N}^{\prime }=\left\{ \xi :\xi \in 
\mathbb{R}^{n}\text{, }\left\vert \xi \right\vert \geq N\right\} . 
\]

From $\left( 2.6\right) $ we deduced that

\[
\left\Vert A^{\alpha }\ast u\right\Vert _{X_{\infty }}\lesssim \left\Vert 
\mathbb{F}^{-1}C\left( \xi ,t\right) \hat{A}^{\alpha }\hat{\varphi}\left(
\xi \right) \right\Vert _{L^{\infty }\left( \Pi _{N}\right) }+ 
\]%
\begin{equation}
\left\Vert \mathbb{F}^{-1}S\left( \xi ,t\right) \hat{A}^{\alpha }\hat{\psi}%
\left( \xi \right) \right\Vert _{L^{\infty }\left( \Pi _{N}\right)
}+\left\Vert \mathbb{F}^{-1}C\left( \xi ,t\right) \hat{A}^{\alpha }\hat{%
\varphi}\left( \xi \right) \right\Vert _{L^{\infty }\left( \Pi _{N}^{\prime
}\right) }+  \tag{2.8}
\end{equation}%
\[
\left\Vert \mathbb{F}^{-1}S\left( \xi ,t\right) \hat{A}^{\alpha }\hat{\psi}%
\left( \xi \right) \right\Vert _{L^{\infty }\left( \Pi _{N}^{\prime }\right)
}+\frac{1}{2}\left\Vert \mathbb{F}^{-1}\hat{A}^{\alpha }\tilde{Q}\left( \xi
,t\right) \hat{g}\left( \xi ,\tau \right) \right\Vert _{L^{\infty }\left(
\Pi _{N}\right) }+ 
\]%
\[
\frac{1}{2}\left\Vert \mathbb{F}^{-1}\hat{A}^{\alpha }\tilde{Q}\left( \xi
,t\right) \hat{g}\left( \xi ,\tau \right) \right\Vert _{L^{\infty }\left(
\Pi _{N}^{\prime }\right) }. 
\]

\bigskip By virtue of Remaks 2.1, 2.2 and properties of sectorial operators
we have the following uniform estimate%
\[
\left\Vert \mathbb{F}^{-1}\hat{A}^{\alpha }\tilde{Q}\left( \xi ,t\right) 
\hat{g}\left( \xi ,\tau \right) \right\Vert _{L^{\infty }\left( \Pi
_{N}\right) }\leq C\left\Vert g\right\Vert _{X_{1}}. 
\]

Hence, due to uniform boundedness of operator functions $C\left( \xi
,t\right) $, $S\left( \xi ,t\right) $, in view of $\left( 2.3\right) $ and\
by Minkowski's inequality for integrals\ we get the uniform estimate 
\[
\left\Vert \mathbb{F}^{-1}C\left( \xi ,t\right) \hat{A}^{\alpha }\hat{\varphi%
}\left( \xi \right) \right\Vert _{L^{\infty }\left( \Pi _{N}\right)
}+\left\Vert \mathbb{F}^{-1}S\left( \xi ,t\right) \hat{A}^{\alpha }\hat{\psi}%
\left( \xi \right) \right\Vert _{L^{\infty }\left( \Pi _{N}\right) }\lesssim 
\]%
\[
\left[ \left\Vert A^{\alpha }\varphi \right\Vert _{X_{1}}+\left\Vert
A^{\alpha }\psi \right\Vert _{X_{1}}+\left\Vert g\right\Vert _{X_{1}}\right]
.
\]%
Moreover, from $\left( 2.6\right) $ we deduced that 
\[
\left\Vert \mathbb{F}^{-1}C\left( \xi ,t\right) \hat{A}^{\alpha }\hat{\varphi%
}\left( \xi \right) \right\Vert _{L^{\infty }\left( \Pi _{N}^{\prime
}\right) }+\left\Vert \mathbb{F}^{-1}S\left( \xi ,t\right) \hat{A}^{\alpha }%
\hat{\psi}\left( \xi \right) \right\Vert _{L^{\infty }}\lesssim 
\]%
\[
\left\Vert \mathbb{F}^{-1}C\left( \xi ,t\right) \hat{A}^{\alpha }\hat{\varphi%
}\left( \xi \right) \right\Vert _{L^{\infty }}+\left\Vert \mathbb{F}%
^{-1}S\left( \xi ,t\right) \hat{A}^{\alpha }\hat{\psi}\left( \xi \right)
\right\Vert _{L^{\infty }}+
\]%
\[
\left\Vert \mathbb{F}^{-1}S\left( \xi ,t\right) \hat{A}^{\alpha }\tilde{Q}%
\left( \xi ,t\right) \hat{g}\left( \xi ,\tau \right) \right\Vert _{L^{\infty
}}\lesssim 
\]%
\begin{equation}
\left\Vert \mathbb{F}^{-1}\left( 1+\left\vert \xi \right\vert ^{2}\right) ^{-%
\frac{s}{2}}C\left( \xi ,t\right) \left( 1+\left\vert \xi \right\vert
^{2}\right) ^{\frac{s}{2}}\hat{A}^{\alpha }\hat{\varphi}\left( \xi \right)
\right\Vert _{L^{\infty }}+  \tag{2.10}
\end{equation}%
\[
\left\Vert \mathbb{F}^{-1}\left( 1+\left\vert \xi \right\vert ^{2}\right) ^{-%
\frac{s}{2}}S\left( \xi ,t\right) \left( 1+\left\vert \xi \right\vert
^{2}\right) ^{\frac{s}{2}}\hat{A}^{\alpha }\hat{\psi}\left( \xi \right)
\right\Vert _{L^{\infty }}+
\]%
\[
\left\Vert \mathbb{F}^{-1}\left( 1+\left\vert \xi \right\vert ^{2}\right) ^{-%
\frac{s}{2}}S\left( \xi ,t\right) \left( 1+\left\vert \xi \right\vert
^{2}\right) ^{\frac{s}{2}}\hat{A}^{\alpha }\tilde{Q}\left( \xi ,t\right) 
\hat{g}\left( \xi ,\tau \right) \right\Vert _{L^{\infty }},
\]%
here, the space $L^{\infty }\left( \Omega ;E\right) $ is denoted by $%
L^{\infty }$. It is clear to see that

\[
\frac{\partial }{\partial \xi _{k}}\left[ \left( 1+\left\vert \xi
\right\vert ^{2}\right) ^{-\frac{s}{2}}\hat{A}^{\alpha }\left( \xi \right)
C\left( \xi ,t\right) \Phi _{0}\left( \xi \right) \right] =
\]%
\[
-s\xi _{k}\left( 1+\left\vert \xi \right\vert ^{2}\right) ^{-\frac{s}{2}-1}%
\hat{A}^{\alpha }\left( \xi \right) C\left( \xi ,t\right) \Phi _{0}\left(
\xi \right) +
\]%
\[
\left( 1+\left\vert \xi \right\vert ^{2}\right) ^{-\frac{s}{2}}\left[ \frac{%
\partial }{\partial \xi _{k}}\left[ \hat{A}^{\alpha }\left( \xi \right)
C\left( \xi ,t\right) \right] \Phi _{0}\left( \xi \right) +\hat{A}^{\alpha
}\left( \xi \right) C\left( \xi ,t\right) \frac{\partial }{\partial \xi _{k}}%
\Phi _{0}\left( \xi \right) \right] =
\]%
\[
-s\xi _{k}\left( 1+\left\vert \xi \right\vert ^{2}\right) ^{-\frac{s}{2}-1}%
\hat{A}^{\alpha }\left( \xi \right) C\left( \xi ,t\right) \Phi _{0}\left(
\xi \right) +
\]%
\[
\left( 1+\left\vert \xi \right\vert ^{2}\right) ^{-\frac{s}{2}}\left\{ \left[
\frac{it}{4}\hat{A}^{\alpha }\left( \xi \right) \eta _{-}\left( \xi \right)
\left( 2\xi _{k}\hat{a}\left( \xi \right) +\left\vert \xi \right\vert ^{2}%
\frac{\partial }{\partial \xi _{k}}\hat{a}\left( \xi \right) +\frac{\partial 
}{\partial \xi _{k}}\hat{A}\left( \xi \right) \right) +\right. \right. 
\]%
\[
\left. \alpha C\left( \xi ,t\right) \hat{A}^{\alpha -1}\left( \xi \right) 
\frac{\partial }{\partial \xi _{k}}\hat{A}\left( \xi \right) \right] \Phi
_{0}\left( \xi \right) +\left. \hat{A}^{\alpha }\left( \xi \right) C\left(
\xi ,t\right) \frac{\partial }{\partial \xi _{k}}\Phi _{0}\left( \xi \right)
\right\} ,
\]

\begin{equation}
\frac{\partial }{\partial \xi _{k}}\left[ \left( 1+\left\vert \xi
\right\vert ^{2}\right) ^{-\frac{s}{2}}\hat{A}^{\alpha }\left( \xi \right)
S\left( \xi ,t\right) \Phi _{1}\left( \xi \right) \right] =  \tag{2.11}
\end{equation}%
\[
-s\xi _{k}\left( 1+\left\vert \xi \right\vert ^{2}\right) ^{-\frac{s}{2}-1}%
\hat{A}^{\alpha }\left( \xi \right) S\left( \xi ,t\right) \Phi _{1}\left(
\xi \right) +
\]%
\[
\left( 1+\left\vert \xi \right\vert ^{2}\right) ^{-\frac{s}{2}}\left\{ \left[
\frac{t}{4}\hat{A}^{\alpha }\left( \xi \right) \eta _{+}\left( \xi \right)
\left( 2\xi _{k}\hat{a}\left( \xi \right) +\left\vert \xi \right\vert ^{2}%
\frac{\partial }{\partial \xi _{k}}\hat{a}\left( \xi \right) +\frac{\partial 
}{\partial \xi _{k}}\hat{A}\left( \xi \right) \right) +\right. \right. 
\]%
\[
\frac{t}{4i}\hat{A}^{\alpha }\left( \xi \right) \eta _{-}\left( \xi \right)
\left( 2\xi _{k}\hat{a}\left( \xi \right) +\left\vert \xi \right\vert ^{2}%
\frac{\partial }{\partial \xi _{k}}\hat{a}\left( \xi \right) +\frac{\partial 
}{\partial \xi _{k}}\hat{A}\left( \xi \right) \right) \eta ^{-2}\left( \xi
\right) +
\]%
\[
\left. \left. \alpha S\left( \xi ,t\right) \hat{A}^{\alpha -1}\left( \xi
\right) \frac{\partial }{\partial \xi _{k}}\hat{A}\left( \xi \right) \right]
\Phi _{1}\left( \xi \right) +\hat{A}^{\alpha }\left( \xi \right) C\left( \xi
,t\right) \frac{\partial }{\partial \xi _{k}}\Phi _{1}\left( \xi \right)
\right\} ,
\]%
where 
\[
\Phi _{0}\left( \xi \right) =\left[ \hat{A}^{1-\frac{1}{2p}-\varepsilon
_{0}}+\left( 1+\left\vert \xi \right\vert ^{2}\right) ^{s\left( 1-\frac{1}{2p%
}\right) }\right] ^{-1}\text{, }0<\varepsilon _{0}<1-\frac{1}{2p},
\]%
\[
\Phi _{1}\left( \xi \right) =\left[ \hat{A}^{\frac{1}{2}-\frac{1}{2p}%
-\varepsilon }+\left( 1+\left\vert \xi \right\vert ^{2}\right) ^{s\left( 
\frac{1}{2}-\frac{1}{2p}\right) }\right] ^{-1}\text{, }0<\varepsilon _{1}<%
\frac{1}{2}-\frac{1}{2p}.\text{ }
\]%
By assumption on $\hat{A}^{\alpha }\left( \xi \right) $, we have the uniform
estmates 
\[
\left\Vert \hat{A}^{\alpha }\left( \xi \right) C\left( \xi ,t\right) \Phi
_{0}\left( \xi \right) \right\Vert _{B\left( E\right) }\leq C\left\Vert \hat{%
A}^{\alpha }\left( \xi \right) \hat{A}^{-\left( 1-\frac{1}{2p}-\varepsilon
_{0}\right) }\left( \xi \right) \right\Vert _{B\left( E\right) }\leq C_{0},
\]%
\[
\left\Vert \hat{A}^{\frac{1}{2}}\left( \xi \right) \eta ^{-1}\left( \xi
\right) \right\Vert _{B\left( E\right) }\left\Vert \hat{A}^{\alpha }\left(
\xi \right) \hat{A}^{-\frac{1}{2}}\left( \xi \right) \Phi _{1}\left( \xi
\right) \right\Vert _{B\left( E\right) }\leq 
\]%
\[
\left\Vert \hat{A}^{\alpha }\left( \xi \right) S\left( \xi ,t\right) \Phi
_{1}\left( \xi \right) \right\Vert _{B\left( E\right) }\leq C\left\Vert \hat{%
A}^{\alpha }\left( \xi \right) \hat{A}^{-\left( 1-\frac{1}{2p}-\varepsilon
_{0}\right) }\left( \xi \right) \right\Vert _{B\left( E\right) }\leq C_{1}.
\]%
Then by calculating $\frac{\partial }{\partial \xi _{k}}\Phi _{0}\left( \xi
\right) $, $\frac{\partial }{\partial \xi _{k}}\Phi _{1}\left( \xi \right) $%
\ and in view of the assumptions on $\frac{\partial }{\partial \xi _{k}}\hat{%
A}\left( \xi \right) $ we obtain 
\[
\text{ }\hat{A}^{\alpha }\left( \xi \right) \frac{\partial }{\partial \xi
_{k}}\Phi _{0}\left( \xi \right) \in B\left( E\right) \text{, }\hat{A}%
^{\alpha }\left( \xi \right) \frac{\partial }{\partial \xi _{k}}\Phi
_{1}\left( \xi \right) \in B\left( E\right) .
\]%
By assumption (4), in view of $s>1+\frac{n}{p}$\ from $\left( 2.3\right) $, $%
\left( 2.11\right) $ for $\beta =\left( \beta _{1},\beta _{2},...,\beta
_{n}\right) $ and $\beta _{k}\in \left\{ 0,1\right\} $ we have the following
uniform estimates%
\[
\sup\limits_{\xi \in \mathbb{R}^{n},t\in \left[ 0,T\right] }\left\vert \xi
\right\vert ^{\left\vert \beta \right\vert +\frac{n}{p}}\left\Vert D^{\beta }%
\left[ \left( 1+\left\vert \xi \right\vert ^{2}\right) ^{-\frac{s}{2}}\hat{A}%
^{\alpha }\left( \xi \right) C\left( \xi ,t\right) \Phi _{0}\left( \xi
\right) \right] \right\Vert _{B\left( E\right) }\leq C_{1},
\]%
\ 
\begin{equation}
\sup\limits_{\xi \in \mathbb{R}^{n},t\in \left[ 0,T\right] }\left\vert \xi
\right\vert ^{\left\vert \beta \right\vert +\frac{n}{p}}\left\Vert D^{\beta }%
\left[ \left( 1+\left\vert \xi \right\vert ^{2}\right) ^{-\frac{s}{2}}\hat{A}%
^{\alpha }\left( \xi \right) S\left( \xi ,t\right) \Phi _{1}\left( \xi
\right) \right] \right\Vert _{B\left( E\right) }\leq C_{2}  \tag{2.12}
\end{equation}%
Moreover, in view of $\left( 2.12\right) $\ we show that the operator
functions 
\[
\left\vert \xi \right\vert ^{\left\vert \beta \right\vert +\frac{n}{p}%
}D^{\beta }\left[ \left( 1+\left\vert \xi \right\vert ^{2}\right) ^{-\frac{s%
}{2}}\hat{A}^{\alpha }\left( \xi \right) C\left( \xi ,t\right) \Phi
_{0}\left( \xi \right) \right] \text{, }
\]%
\[
\left\vert \xi \right\vert ^{\left\vert \beta \right\vert +\frac{n}{p}%
}D^{\beta }\left[ \left( 1+\left\vert \xi \right\vert ^{2}\right) ^{-\frac{s%
}{2}}\hat{A}^{\alpha }\left( \xi \right) S\left( \xi ,t\right) \Phi
_{1}\left( \xi \right) \right] 
\]%
are uniformly $R$-bounded in $E$. Hence, by Fourier multiplier theorems (see
e.g. $\left[ \text{13, Theorem 4.3}\right] $) we get that the functions $%
\left( 1+\left\vert \xi \right\vert ^{2}\right) ^{-\frac{s}{2}}\hat{A}%
^{\alpha }\left( \xi \right) C\left( \xi ,t\right) \Phi _{i}\left( \xi
\right) $ are $L^{p}\left( \mathbb{R}^{n};E\right) \rightarrow L^{\infty
}\left( \mathbb{R}^{n};E\right) $ Fourier multipliers. Then by Minkowski's
inequality for integrals, from $\left( 2.3\right) $, $\left( 2.10\right) $
and $\left( 2.11\right) -\left( 2.12\right) $ we have 
\[
\left\Vert F^{-1}C\left( \xi ,t\right) \hat{A}^{\alpha }\hat{\varphi}\left(
\xi \right) \right\Vert _{L^{\infty }}+\left\Vert \mathbb{F}^{-1}S\left( \xi
,t\right) \hat{A}^{\alpha }\hat{\psi}\left( \xi \right) \right\Vert
_{L^{\infty }}\lesssim 
\]%
\[
\left\Vert F^{-1}C\left( \xi ,t\right) \eta ^{-2}\hat{\varphi}\right\Vert
_{L^{\infty }}+\left\Vert \mathbb{F}^{-1}S\left( \xi ,t\right) \eta ^{-1}%
\hat{\psi}\right\Vert _{L^{\infty }}\lesssim 
\]

\begin{equation}
\left[ \left\Vert \varphi \right\Vert _{\mathbb{E}_{0p}}+\left\Vert \psi
\right\Vert _{\mathbb{E}_{1p}}+\left\Vert g\right\Vert _{W^{s,p}}\right] . 
\tag{2.13}
\end{equation}%
Moreover, by virtue of Remaks 2.1, 2.2 and by reasoning as the above, we
have the followin estimate 
\begin{equation}
\left\Vert F^{-1}\hat{A}^{\alpha }\tilde{Q}\left( \xi ,t\right) \right\Vert
_{X_{\infty }}\leq C\dint\limits_{0}^{t}\left( \left\Vert g\left( .,\tau
\right) \right\Vert _{W^{s,p}}+\left\Vert g\left( .,\tau \right) \right\Vert
_{X_{1}}\right) d\tau  \tag{2.14}
\end{equation}%
uniformly in $t\in \left[ 0,T\right] $. Thus, from $\left( 2.6\right) $, $%
\left( 2.13\right) $ and $\left( 2.14\right) $ we obtain 
\begin{equation}
\left\Vert A^{\alpha }\ast u\right\Vert _{X_{\infty }}\leq C\left[
\left\Vert \varphi \right\Vert _{\mathbb{E}_{0p}}+\left\Vert A^{\alpha
}\varphi \right\Vert _{X_{1}}\right. +  \tag{2.15}
\end{equation}

\[
\left\Vert \psi \right\Vert _{\mathbb{E}_{1p}}+\left\Vert A^{\alpha }\psi
\right\Vert _{X_{1}}+\left. \dint\limits_{0}^{t}\left( \left\Vert g\left(
.,\tau \right) \right\Vert _{Y^{s,p}}+\left\Vert g\left( .,\tau \right)
\right\Vert _{X_{1}}\right) d\tau \right] . 
\]%
By differentiating $\left( 2.6\right) $, in a similar way we get 
\begin{equation}
\left\Vert A^{\alpha }\ast u_{t}\right\Vert _{X_{\infty }}\leq C\left[
\left\Vert \varphi \right\Vert _{\mathbb{E}_{0p}}+\left\Vert A^{\alpha }\ast
\varphi \right\Vert _{X_{1}}\right. +  \tag{2.16}
\end{equation}

\[
\left\Vert A^{\alpha }\ast \psi \right\Vert _{\mathbb{E}_{1p}}+\left\Vert
A^{\alpha }\ast \psi \right\Vert _{X_{1}}+\left. \dint\limits_{0}^{t}\left(
\left\Vert g\left( .,\tau \right) \right\Vert _{Y^{s,p}}+\left\Vert g\left(
.,\tau \right) \right\Vert _{X_{1}}\right) d\tau \right] . 
\]

Then from $\left( 2.15\right) $ and $\left( 2.16\right) $ in view of Remarks
2.1, 2.2 we obtain the assertion.

\textbf{Theorem 2.2. }Let the Condition 2.1 holds, $s>1+\frac{n}{p}$ and let 
$0<\alpha <1-\frac{1}{2p}$. Then for $\varphi \in $ $\mathbb{E}_{0p}$, $\psi
\in $ $\mathbb{E}_{1p}$ and $g\in Y^{s,p}$ the problem $\left( 2.1\right)
-\left( 2.2\right) $ has a unique generalized solution $u\in C^{2}\left( %
\left[ 0,T\right] ;Y^{s,p}\right) $ and the following uniform estimate holds%
\begin{equation}
\left( \left\Vert A^{\alpha }\ast u\right\Vert _{Y^{s,p}}+\left\Vert
A^{\alpha }\ast u_{t}\right\Vert _{Y^{s,p}}\right) \leq   \tag{2.17}
\end{equation}

\[
C_{0}\left[ \left\Vert \varphi \right\Vert _{\mathbb{E}_{0p}}+\left\Vert
\psi \right\Vert _{\mathbb{E}_{1p}}+\dint\limits_{0}^{t}\left\Vert g\left(
.,\tau \right) \right\Vert _{Y^{s,p}}d\tau \right] . 
\]%
\textbf{Proof. }From $\left( 2.5\right) $ and $\left( 2.11\right) $ we get
the following uniform estimate 
\begin{equation}
\left( \left\Vert \mathbb{F}^{-1}\left( 1+\left\vert \xi \right\vert
^{2}\right) ^{\frac{s}{2}}\hat{A}^{\alpha }\hat{u}\right\Vert
_{X_{p}}+\left\Vert \mathbb{F}^{-1}\left( 1+\left\vert \xi \right\vert
^{2}\right) ^{\frac{s}{2}}\hat{A}^{\alpha }\hat{u}_{t}\right\Vert
_{X_{p}}\right) \leq  \tag{2.18}
\end{equation}

\[
C\left\{ \left\Vert \mathbb{F}^{-1}\left( 1+\left\vert \xi \right\vert
^{2}\right) ^{\frac{s}{2}}C\left( \xi ,t\right) \hat{A}^{\alpha }\hat{\varphi%
}\right\Vert _{X_{p}}\right. +\left\Vert \mathbb{F}^{-1}\left( 1+\left\vert
\xi \right\vert ^{2}\right) ^{\frac{s}{2}}\hat{A}^{\alpha }S\left( \xi
,t\right) \hat{\psi}\right\Vert _{X_{p}}+ 
\]

\[
\left. \dint\limits_{0}^{t}\left\Vert \left( 1+\left\vert \xi \right\vert
^{2}\right) ^{\frac{s}{2}}\hat{A}^{\alpha }\tilde{Q}\left( \xi ,t\right) 
\hat{g}\left( \xi ,\tau \right) \right\Vert _{X_{p}}d\tau \right\} . 
\]

\bigskip By using the Fourier multiplier theorem $\left[ \text{13, Theorem
4.3}\right] $ and by reasoning as in Theorem 2.1\ we get that $C\left( \xi
,t\right) $, $S\left( \xi ,t\right) $ and $\hat{A}^{\alpha }S\left( \xi
,t\right) $\ are Fourier multipliers in $L^{p}\left( \mathbb{R}^{n};E\right) 
$ uniformly with respect to $t\in \left[ 0,T\right] $. So, the estimate $%
\left( 2.18\right) $ by using the Minkowski's inequality for integrals
implies $\left( 2.17\right) .$

\begin{center}
\textbf{3. Local well posedness of IVP for nonlinear nonlocal WE}
\end{center}

In this section, we will show the local existence and uniqueness of solution
for the nonlinear problem $(1.1)-(1.2)$.

For this aim we need the following lemmas. Here, we will denote $L^{p}\left( 
\mathbb{R}^{n};E\right) $, $W^{s,p}\left( \mathbb{R}^{n};E\right) $ by by $%
X_{p}$ and $Y^{s,p}$, respectively. Here, we assume that $E$ is a Banach
algebra. By reasoning as in $[$8, 13, 26$]$, we show the following lemmas
concerning the behaviour of the nonlinear term in $E-$valued space $Y^{s,p}$.

\textbf{\ Lemma 3.1.} Let $s\geq 0$, $f\in C^{\left[ s\right] +1}\left( 
\mathbb{R};E\right) $ with $f(0)=0$. Then for any $u\in Y^{s,p}\cap
L^{\infty }$, we have $f(u)\in Y^{s,p}\cap X_{\infty }$. Moreover, there is
some constant $A(M)$ depending on $M$ such that for all $u\in Y^{s,p}\cap
L^{\infty }$ with $\left\Vert u\right\Vert _{X_{\infty }}\leq M,$%
\begin{equation}
\left\Vert f(u)\right\Vert _{Y^{s,p}}\leq C\left( M\right) \left\Vert
u)\right\Vert _{Y^{s,p}}.  \tag{3.1}
\end{equation}%
\textbf{Proof. }For $s=0$ in view of $f(0)=0$, we get%
\[
f\left( u\right) =u\dint\limits_{0}^{1}f\left( \sigma u\right) d\sigma . 
\]

It follows that 
\[
\left\Vert f\left( u\right) \right\Vert _{X_{p}}\leq C\left( M\right)
\left\Vert u\right\Vert _{X_{p}}. 
\]

If $s>0$ is a positive integer, we have 
\begin{equation}
\left\Vert f(u)\right\Vert _{Y^{s,p}}\leq C\left[ \left\Vert f(u)\right\Vert
_{X_{p}}+\dsum\limits_{k=1}^{n}\left\Vert \frac{\partial ^{s}}{\partial x_{k}%
}f(u)\right\Vert _{X_{p}}\right] .  \tag{3.2}
\end{equation}%
By calculation of derivativie and applying Holder inequality we get%
\[
\left\Vert \frac{\partial ^{s}}{\partial x_{i}}f(u)\right\Vert _{X_{p}}\leq
\dsum\limits_{l=1}^{s}\dsum\limits_{\alpha }\left\Vert f^{\left( l\right)
}(u)\frac{\partial ^{\beta _{1}}u}{\partial x_{i}}\frac{\partial ^{\beta
_{2}}u}{\partial x_{i}}...\frac{\partial ^{\beta _{l}}u}{\partial x_{i}}%
\right\Vert _{X_{p}}\leq 
\]%
\begin{equation}
\dsum\limits_{l=1}^{s}\dsum\limits_{\alpha }\left\Vert f^{\left( l\right)
}(u)\right\Vert _{X_{\infty }}\dprod\limits_{k=1}^{l}\left\Vert \frac{%
\partial ^{\beta _{k}}u}{\partial x_{i}}\right\Vert _{X_{p_{k}}}\text{, }%
i=1,2,...,n,  \tag{3.3}
\end{equation}%
where 
\[
\beta =\left( \beta _{1},\beta _{2},,,,.\beta _{l}\right) \text{, }\beta
_{k}\geq 1\text{, }\beta _{1}+\beta _{2}+...+\beta _{l}=l\text{, }p_{k}=%
\frac{pl}{\beta _{k}}. 
\]

Applying Gagliardo-Nirenberg's inequality in $E$-valued $X_{p}$ spaces, we
have 
\begin{equation}
\left\Vert \frac{\partial ^{\beta _{k}}u}{\partial x_{i}}\right\Vert
_{X_{p_{k}}}\leq C\left\Vert u\right\Vert _{X_{\infty }}^{1-\frac{\beta _{k}%
}{l}}\left\Vert \frac{\partial ^{s}u}{\partial x_{i}^{s}}\right\Vert
_{X_{p}}^{\frac{\beta _{k}}{l}}.  \tag{3.4}
\end{equation}

Hence, from $\left( 3.3\right) $ and $\left( 3.4\right) $ we deduced%
\begin{equation}
\left\Vert \frac{\partial ^{s}}{\partial x_{i}}f(u)\right\Vert _{X_{p}}\leq
C\left( M\right) \left\Vert \frac{\partial ^{s}u}{\partial x_{i}^{s}}%
\right\Vert _{X_{p}}.  \tag{3.5}
\end{equation}

Then combinig $\left( 3.2\right) $, $\left( 3.3\right) $ and $\left(
3.5\right) $ we obtain $\left( 3.1\right) $.

If $s$ is not integer number, let $m=\left[ s\right] $. From the above
proof, we have 
\[
\left\Vert f(u)\right\Vert _{Y^{m,p}}\leq C\left( M\right) \left\Vert
u)\right\Vert _{Y^{m,p}}\text{, }\left\Vert f(u)\right\Vert _{Y^{m+1,p}}\leq
C\left( M\right) \left\Vert u)\right\Vert _{Y^{m+1,p}}. 
\]

Then using interpolation between $W^{m+1,p}$ and $W^{m,p}$ yields $\left(
3.1\right) $ for all $s\geq 0.$

By using Lemma 3.1 and properties of convolution operators we obtain

\textbf{Corollary 3.1. }Let $s\geq 0$, $f\in C^{\left[ s\right] +1}\left( 
\mathbb{R};E\right) $ with $f(0)=0$. Moreover, assume $\Phi \in L^{\infty
}\left( \mathbb{R}^{n};B\left( E\right) \right) $. Then for any $u\in
Y^{s,p}\cap L^{\infty }$, we have $f(u)\in Y^{s,p}\cap X_{\infty }$.
Moreover, there is some constant $A(M)$ depending on $M$ such that for all $%
u\in Y^{s,p}\cap L^{\infty }$ with $\left\Vert u\right\Vert _{X_{\infty
}}\leq M,$%
\[
\left\Vert \Phi \ast f(u)\right\Vert _{Y^{s,p}}\leq C\left( M\right)
\left\Vert u)\right\Vert _{Y^{s,p}}. 
\]

\textbf{Lemma 3.2. }Let $s\geq 0,$ $f\in C^{\left[ s\right] +1}\left( 
\mathbb{R};E\right) $. Then for for any $M$ there is some constant $K(M)$
depending on $M$ such that for all $u$, $\upsilon \in Y^{s,p}\cap X_{\infty
} $ with $\left\Vert u\right\Vert _{X_{\infty }}\leq M$, $\left\Vert
\upsilon \right\Vert _{X_{\infty }}\leq M$, $\left\Vert u\right\Vert
_{Y^{s,p}}\leq M$, $\left\Vert \upsilon \right\Vert _{Y^{s,p}}\leq M,$%
\[
\left\Vert f(u)-f(\upsilon \right\Vert _{Y^{s,p}}\leq K\left( M\right)
\left\Vert u-\upsilon \right\Vert _{Y^{s,p}},\text{ }\left\Vert
f(u)-f(\upsilon \right\Vert _{X_{\infty }}\leq K\left( M\right) \left\Vert
u-\upsilon \right\Vert _{X_{\infty }}. 
\]

By reasoning as in $\left[ \text{13, Lemma 3.4}\right] $ and $\left[ \text{%
26, Lemma X 4}\right] $ we have, respectively

\textbf{Corollary 3.2.} Let $s>\frac{n}{2}$, $f\in C^{\left[ s\right]
+1}\left( \mathbb{R};E\right) $. Then for any $M$ there is a constant $K(M)$
depending on $M$ such that for all $u$, $\upsilon \in Y^{s,p}$ with $%
\left\Vert u\right\Vert _{Y^{s,p}}\leq M$, $\left\Vert \upsilon \right\Vert
_{Y^{s,p}}\leq M,$%
\[
\left\Vert f(u)-f(\upsilon \right\Vert _{Y^{s,p}}\leq K\left( M\right)
\left\Vert u-\upsilon \right\Vert _{Y^{s,p}}. 
\]

\bigskip \textbf{Lemma 3.3. }If $s>0$, then $Y_{\infty }^{s,p}$ is an
algebra. Moreover, for \ $f$, $g\in Y_{\infty }^{s,p},$ 
\[
\left\Vert fg\right\Vert _{Y^{s,p}}\leq C\left[ \left\Vert f\right\Vert
_{X_{\infty }}+\left\Vert g\right\Vert _{Y^{s,p}}+\left\Vert f\right\Vert
_{Y^{s,p}}+\left\Vert g\right\Vert _{X_{\infty }}\right] . 
\]

\textbf{\ }By using, The Corollary 3.1 and Lemma 3.3 we obta\i n

\textbf{Lemma 3.4.} Let $s\geq 0$, $f\in C^{\left[ s\right] +1}\left( 
\mathbb{R};E\right) $ and $f\left( u\right) =O\left( \left\vert u\right\vert
^{\gamma +1}\right) $ for $u\rightarrow 0$, $\gamma \geq 1$ be a positive
integer. If $u\in Y_{\infty }^{s,p}$ and $\left\Vert u\right\Vert
_{X_{\infty }}\leq M$, then 
\[
\left\Vert f(u)\right\Vert _{Y^{s,p}}\leq C\left( M\right) \left[ \left\Vert
u\right\Vert _{Y^{s,p}}\left\Vert u\right\Vert _{X_{\infty }}^{\gamma }%
\right] , 
\]%
\[
\left\Vert f(u)\right\Vert _{X_{1}}\leq C\left( M\right) \left\Vert
u\right\Vert _{X_{p}}^{p}\left\Vert u\right\Vert _{X_{\infty }}^{\gamma -1}. 
\]

\textbf{Corollary 3.3. }Let $s\geq 0$, $f\in C^{\left[ s\right] +1}\left( 
\mathbb{R};E\right) $ and $f\left( u\right) =O\left( \left\vert u\right\vert
^{\gamma +1}\right) $ for $u\rightarrow 0$, $\gamma \geq 1$ be a positive
integer. Moreover, assume $\Phi \in L^{\infty }\left( \mathbb{R}^{n};B\left(
E\right) \right) .$ If $u\in Y_{\infty }^{s,p}$ and $\left\Vert u\right\Vert
_{X_{\infty }}\leq M$, then 
\[
\left\Vert \Phi \ast f(u)\right\Vert _{Y^{s,p}}\leq C\left( M\right) \left[
\left\Vert u\right\Vert _{Y^{s,p}}\left\Vert u\right\Vert _{X_{\infty
}}^{\gamma }\right] , 
\]%
\[
\left\Vert \Phi \ast f(u)\right\Vert _{X_{1}}\leq C\left( M\right)
\left\Vert u\right\Vert _{X_{p}}^{p}\left\Vert u\right\Vert _{X_{\infty
}}^{\gamma -1}. 
\]

\textbf{Lemma 3.5.} Let $s\geq 0$, $f\in C^{\left[ s\right] +1}\left( 
\mathbb{R};E\right) $ and $f\left( u\right) =O\left( \left\vert u\right\vert
^{\gamma +1}\right) $\ for $u\rightarrow 0$, $\gamma \geq 0$ be a positive
integer. If $u,$ $\upsilon \in Y_{\infty }^{s,p}$, $\left\Vert u\right\Vert
_{Y^{s,p}}\leq M$, $\left\Vert \upsilon \right\Vert _{Y^{s,p}}\leq M$ and $%
\left\Vert u\right\Vert _{X_{\infty }}\leq M$, $\left\Vert \upsilon
\right\Vert _{X_{\infty }}\leq M$, then 
\[
\left\Vert f(u)-f(\upsilon )\right\Vert _{Y^{s,p}}\leq C\left( M\right) 
\left[ \left( \left\Vert u\right\Vert _{X_{\infty }}-\left\Vert \upsilon
\right\Vert _{X_{\infty }}\right) \left( \left\Vert u\right\Vert
_{Y^{s,p}}+\left\Vert \upsilon \right\Vert _{Y^{s,p}}\right) \right. 
\]%
\[
\left( \left\Vert u\right\Vert _{X_{\infty }}+\left\Vert \upsilon
\right\Vert _{X_{\infty }}\right) ^{\gamma -1}, 
\]%
\[
\left\Vert f(u)-f(\upsilon \right\Vert _{X_{1}}\leq C\left( M\right) \left(
\left\Vert u\right\Vert _{X_{\infty }}+\left\Vert \upsilon \right\Vert
_{X_{\infty }}\right) ^{\gamma -1}\left( \left\Vert u\right\Vert
_{X_{p}}+\left\Vert \upsilon \right\Vert _{X_{p}}\right) \left\Vert
u-\upsilon \right\Vert _{X_{p}}. 
\]

Let $E_{0}$ denotes the real interpolation space between $Y^{s,p}\left(
A,E\right) $ and $X_{p}$ with $\theta =\frac{1}{2p}$, i.e. 
\[
\text{ }\mathbb{E}_{0p}=\left( Y^{s,p}\left( A,E\right) ,X_{p}\right) _{%
\frac{1}{2p},p}. 
\]

Here, $Y_{0}\left( A^{\alpha }\right) $ and $Y_{1}\left( A^{\alpha }\right) $
are the spaces defined in Section 2.

\textbf{Remark 3.1. }By using J.Lions-I. Petree result (see e.g $\left[ 
\text{27, \S\ 1.8}\right] $) we obtain that the map $u\rightarrow u\left(
t_{0}\right) $, $t_{0}\in \left[ 0,T\right] $ is continuous and surjective
from $Y^{2,s,p}\left( A,E\right) $ onto $\mathbb{E}_{0p}$ and there is a
constant $C_{1}$ such that 
\begin{equation}
\left\Vert u\left( t_{0}\right) \right\Vert _{\mathbb{E}_{0p}}\leq
C_{1}\left\Vert u\right\Vert _{Y^{2,s,p}\left( A,E\right) }\text{, }1\leq
p\leq \infty \text{.}  \tag{3.6}
\end{equation}

Let 
\[
C^{2}\left( Y_{1}^{s,p}\left( A\right) \right) =C^{\left( 2\right) }\left( 
\left[ 0,T\right] ;Y_{1}^{s,p}\left( A,E\right) \right) \text{, }%
C^{2,s}\left( A,E\right) =C^{\left( 2\right) }\left( \left[ 0,T\right]
;Y^{s,p}\left( A,E\right) \right) . 
\]

\textbf{Definition 3.1. }Let $T>0$ and $\varphi \in $ $Y_{0}\left( A^{\alpha
}\right) $, $\psi \in $ $Y_{1}\left( A^{\alpha }\right) $. The function $u$ $%
\in C^{2}\left( Y_{1}^{s,p}\left( A\right) \right) $ satisfies of the
problem $(1.1)-(1.2)$ is called the continuous solution\ or the strong
solution of $(1.1)-(1.2)$. If $T<\infty $, then $u\left( x,t\right) $ is
called the local strong solution of the problem $(1.1)-(1.2)$. If $T=\infty $%
, then $u\left( x,t\right) $ is called the global strong solution of $%
(1.1)-(1.2)$.

\textbf{Condition 3.1. }Assume:

(1) the Condition 2.1 holds for $s>\frac{n}{p}$ and $0<\alpha <1-\frac{1}{2p}
$;

(2) the kernel $g=g\left( x\right) $ is a bounded integrable operator
function in $E$ such that $\Delta g\in L^{\infty }\left( \mathbb{R}%
^{n};B\left( E\right) \right) $;

(3) the function $u\rightarrow f\left( u\right) $: continuous from $u\in 
\mathbb{E}_{0p}$ into $E$, $f\in C^{k}\left( \mathbb{R};E\right) $ with $k$
an integer, $k\geq s>\frac{n}{p}$ and $f\left( u\right) =O\left( \left\vert
u\right\vert ^{\gamma +1}\right) $ for $u\rightarrow 0$, $\gamma \geq 1$ be
a positive integer.

\textbf{Remark 3.2. }We will use Lemmas 3.1-3.5 and Corollary 3.3 in the
follwing results. Note that, inspite of in these Lemmas and Corollary were
assumed $E$ to be Banach algebra, here it is sufficient to take $E$ UMD
space. Really, since the solution $u$ of $\left( 1.1\right) -\left(
1.2\right) $ is assumed to be from the space $Y^{2,s,p}\left( A,E\right) $.
Then by assumption (3) of the Condition 3.1 and by Remarke 1.1 the function $%
u\rightarrow f\left( u\right) $ is continuous from $u\in Y^{2,s,p}$ into $E$%
. Hence, Lemmas 3.1-3.5 and Corollary 3.3 are yield for $u\in
Y^{2,s,p}\left( A,E\right) $, when $E$ is only UMD spaces.

Let 
\[
\hat{Y}_{1}^{s,p}\left( A^{\alpha };E\right) =\hat{Y}^{s,p}\left( A^{\alpha
};E\right) \cap X_{1}\left( A^{\alpha }\right) \text{, }\hat{Y}^{s,p}\left(
A^{\alpha };E\right) =\left\{ u\in Y^{s,p}\left( A^{\alpha };E\right) \text{,%
}\right. 
\]%
\[
\text{ }\left\Vert u\right\Vert _{\hat{Y}^{s,p}\left( A^{\alpha };E\right)
}=\left\Vert A^{\alpha }\ast u\right\Vert _{X_{p}}+\left. \left\Vert \mathbb{%
F}^{-1}\left( 1+\left\vert \xi \right\vert ^{2}\right) ^{\frac{s}{2}}\hat{u}%
\right\Vert _{X_{p}}<\infty \right\} . 
\]

Main aim of this section is to prove the following results:

\textbf{Theorem 3.1. }Let the Condition 3.1 holds. Then there exists a
constant $\delta >0$ such that for any $\varphi \in $ $Y_{0}\left( A^{\alpha
}\right) $ and $\psi \in $ $Y_{1}\left( A^{\alpha }\right) $ satisfying 
\begin{equation}
\left\Vert \varphi \right\Vert _{\mathbb{E}_{0p}}+\left\Vert A^{\alpha }\ast
\varphi \right\Vert _{X_{1}}+\left\Vert \psi \right\Vert _{\mathbb{E}%
_{1p}}+\left\Vert A^{\alpha }\ast \psi \right\Vert _{X_{1}}\leq \delta , 
\tag{3.7}
\end{equation}%
problem $\left( 1.1\right) -\left( 1.2\right) $ has a unique local strange
solution $u\in C^{2}\left( Y_{1}^{s,p}\left( A\right) \right) $. Moreover,

\begin{equation}
\sup_{t\in \left[ 0,T\right] }\left( \left\Vert u\left( .,t\right)
\right\Vert _{\hat{Y}_{1}^{s,p}\left( A^{\alpha },E\right) }+\left\Vert
u_{t}\left( .,t\right) \right\Vert _{\hat{Y}_{1}^{s,p}\left( A^{\alpha
};E\right) }\right) \leq C\delta ,  \tag{3.8}
\end{equation}%
where the constant $C$ depends only on $A$, $E$, $g$, $f$ and initial values.

\textbf{Proof.} By $(2.5)$, $\left( \left( 2.6\right) \right) $ the problem
of finding a solution $u$ of $(1.1)-\left( 1.2\right) $ is equivalent to
finding a fixed point of the mapping

\begin{equation}
G\left( u\right) =C_{1}\left( t\right) \varphi \left( x\right) +S_{1}\left(
t\right) \psi \left( x\right) +Q\left( u\right) ,  \tag{3.9}
\end{equation}%
where $C_{1}\left( t\right) $, $S_{1}\left( t\right) $ are defined by $%
\left( 2.6\right) $ and $Q\left( u\right) $ is a map defined by 
\[
Q\left( u\right) =-\dint\limits_{0}^{t}\mathbb{F}^{-1}\left[ U\left( \xi
,t-\tau \right) \left\vert \xi \right\vert ^{2}\hat{g}\left( \xi \right) 
\hat{f}\left( u\right) \left( \xi ,\tau \right) \right] d\tau . 
\]%
We define the metric space 
\[
C\left( T,A\right) =C_{\delta }^{2}\left( Y_{1}^{s,p}\left( A\right) \right)
=\left\{ u\in C^{2,s}\left( A,E\right) \text{, }\left\Vert u\right\Vert
_{C^{2,s,p}\left( T,A\right) }\leq 5C_{0}\delta \right\} 
\]%
equipped with the norm defined by 
\[
\left\Vert u\right\Vert _{C\left( T,A\right) }=\sup\limits_{t\in \left[ 0,T%
\right] }\left[ \left\Vert A^{\alpha }\ast u\left( .,t\right) \right\Vert
_{X_{\infty }}+\left\Vert u\left( .,t\right) \right\Vert _{Y^{s,p}}+\right. 
\]%
\[
\left. \left\Vert A^{\alpha }\ast u_{t}\left( .,t\right) \right\Vert
_{X_{\infty }}+\left\Vert u_{t}\left( .,t\right) \right\Vert _{Y^{s,p}} 
\right] , 
\]%
where $\delta >0$ satisfies $\left( 3.7\right) $ and $C_{0}$ is a constant
in Theorem 2.1 and 2.2. It is easy to prove that $C\left( T,A\right) $ is a
complete metric space. From imbedding in Sobolev-Lions space $Y^{s,p}\left(
A,E\right) $ (see e.g. $\left[ \text{30}\right] $, Theorem 1) and trace
result $\left( 3.6\right) $ we got that $\left\Vert u\right\Vert _{X_{\infty
}}\leq 1$ if we take that $\delta $ is enough small. For $\varphi \in $ $%
Y_{0}\left( A^{\alpha }\right) $ and $\psi \in $ $Y_{1}\left( A^{\alpha
}\right) $, let 
\[
\left\Vert \varphi \right\Vert _{\mathbb{E}_{0p}}+\left\Vert A^{\alpha }\ast
\varphi \right\Vert _{X_{1}}+\left\Vert \psi \right\Vert _{\mathbb{E}%
_{1p}}+\left\Vert A^{\alpha }\ast \psi \right\Vert _{X_{1}}=\delta . 
\]%
So, we will find $T$ and $M$ so that $G$ is a contraction on $%
C^{2,s,p}\left( T,A\right) $. By Theorems 2.1, 2.2 and Corollary 3.3 $\Delta
g\ast f\left( u\right) \in Y_{1}^{s,p}$. So, problem $(1.1)-\left(
1.2\right) $ has a solution satisfies the following 
\begin{equation}
G\left( u\right) \left( x,t\right) =C_{1}\left( t\right) \varphi
+S_{1}\left( t\right) \psi +Qu,  \tag{3.10}
\end{equation}%
where $C_{1}\left( t\right) $, $S_{1}\left( t\right) $ are defined by $%
\left( 2.5\right) $ and $\left( 2.6\right) $. By assumptions, it is easy to
see that the map $G$ is well defined for $f\in C^{\left[ s\right] +1}\left( 
\mathbb{E}_{0p};E\right) $. First,\ let us prove that the map $G$ has a
unique fixed point in $C\left( T,A\right) $. For this aim, it is sufficient
to show that the operator $G$ maps $C\left( T,A\right) $ into $C\left(
T,A\right) $ and $G$ is strictly contractive if $\delta $ is suitable
small.\ In fact, by $(2.7)$ in Theorem 2.1, Corollary 3.3 and in view of $%
\left( 3.7\right) $, we have%
\begin{equation}
\left\Vert A^{\alpha }\ast G\left( u\right) \right\Vert _{X_{\infty
}}+\left\Vert A^{\alpha }\ast G\left( u\right) _{t}\right\Vert _{X_{\infty
}}\leq 2C_{0}\left[ \left\Vert \varphi \right\Vert _{Y_{0}^{\alpha }\left(
A^{\alpha }\right) }+\right.  \tag{3.11}
\end{equation}

\[
\left\Vert \psi \right\Vert _{Y_{1}^{\alpha }\left( A^{\alpha }\right)
}+\left. \dint\limits_{0}^{t}\left( \left\Vert \left[ \Delta g\ast f\left(
\left( u\right) \right) \right] \right\Vert _{Y^{s,p}}+\left\Vert \left[
\Delta g\ast f\left( \left( u\right) \right) \right] \right\Vert
_{X_{1}}\right) d\tau \right] \leq 
\]%
\[
2C_{0}\delta +C\dint\limits_{0}^{t}\left( \left\Vert u\left( \tau \right)
\right\Vert _{Y^{s,p}}\left\Vert u\left( \tau \right) \right\Vert
_{X_{\infty }}^{\gamma }+\left\Vert u\left( \tau \right) \right\Vert
_{X_{p}}^{p}\left\Vert u\left( \tau \right) \right\Vert _{X_{\infty
}}^{\gamma -1}\right) d\tau \leq 
\]%
\[
2C_{0}\delta +C\left\Vert u\right\Vert _{C^{2,s,p}\left( T,A\right)
}^{\gamma +1}. 
\]%
On the oher hand, by $(2.17)$ in Theorem 2.2, Corollary 3.3 and $\left(
3.7\right) $, we get

\begin{equation}
\left( \left\Vert A^{\alpha }\ast G\left( u\right) \right\Vert
_{Y^{s,p}}+\left\Vert A^{\alpha }\ast G\left( u\right) _{t}\right\Vert
_{Y^{s,p}}\right) \leq  \tag{3.12}
\end{equation}

\[
2C_{0}\left( \left\Vert \varphi \right\Vert _{\mathbb{E}_{0p}}+\left\Vert
\psi \right\Vert _{\mathbb{E}_{1p}}+\dint\limits_{0}^{t}\left\Vert \Delta %
\left[ g\ast f\left( \left( u\right) \right) \right] \right\Vert
_{Y^{s,p}}d\tau \right) \leq 
\]%
\[
2C_{0}\delta +\dint\limits_{0}^{t}\left[ \left\Vert u\left( \tau \right)
\right\Vert _{Y^{s,p}}\left\Vert u\left( \tau \right) \right\Vert
_{X_{\infty }}^{\gamma }\right] d\tau \leq 2C_{0}\delta +C\left\Vert
u\right\Vert _{C^{2,s,p}\left( T,A\right) }^{\gamma +1}. 
\]

Hence, combining $\left( 3.11\right) $ with $\left( 3.12\right) $ we obtain 
\begin{equation}
\left\Vert A^{\alpha }\ast G\left( u\right) \right\Vert _{Y_{\infty
}^{s,p}}+\left\Vert A^{\alpha }\ast G\left( u\right) _{t}\right\Vert
_{Y_{\infty }^{s,p}}\leq 4C_{0}\delta +C\left\Vert u\right\Vert
_{C^{2,s,p}\left( T,A\right) }^{\gamma +1}.  \tag{3.13}
\end{equation}%
Therefore, taking that $\delta $ is enough small such that $C\left(
5C_{8}\delta \right) ^{\gamma }<\frac{1}{5}$, then by Theorems 2.1, 2.2 and $%
\left( 3.13\right) $, $G$ maps $C\left( T,A\right) $ into $C\left(
T,A\right) $.

Now, we are going to prove that the map $G$ is strictly contractive. Let $%
u_{1}$, $u_{2}\in $ $C\left( T,A\right) $ given. From $\left( 3.10\right) $
we get%
\[
G\left( u_{1}\right) -G\left( u_{2}\right) = 
\]%
\[
\dint\limits_{0}^{t}\left[ S\left( x,t-\tau \right) \Delta g\ast \left(
f\left( u_{1}\right) \left( \tau \right) -f\left( u_{2}\right) \left( \tau
\right) \right) \right] d\tau \text{, }t\in \left( 0,T\right) . 
\]

\bigskip By $(2.7)$ in Theorem 2.1 and Corollary 3.3, we have%
\begin{equation}
\left\Vert A^{\alpha }\ast \left[ G\left( u_{1}\right) -G\left( u_{2}\right) %
\right] \right\Vert _{X_{\infty }}+\left\Vert A^{\alpha }\ast \left[ G\left(
u_{1}\right) -G\left( u_{2}\right) \right] _{t}\right\Vert _{X_{\infty }}\leq
\tag{3.14}
\end{equation}

\[
\dint\limits_{0}^{t}\left( \left\Vert \left[ \Delta g\ast \right] \left[
f\left( u_{1}\right) -f\left( u_{2}\right) \right] \right\Vert
_{Y^{s,p}}+\left\Vert \Delta g\ast \left[ f\left( u_{1}\right) -f\left(
u_{2}\right) \right] \right\Vert _{X_{1}}\right) d\tau \leq 
\]%
\[
\dint\limits_{0}^{t}\left\{ \left\Vert u_{1}-u_{2}\right\Vert _{X_{\infty
}}\left( \left\Vert u_{1}\right\Vert _{Y^{s,p}}+\left\Vert u_{2}\right\Vert
_{Y^{s,p}}\right) \right. \left( \left\Vert u_{1}\right\Vert _{X_{\infty
}}+\left\Vert u_{2}\right\Vert _{X_{\infty }}\right) ^{\gamma -1}+ 
\]%
\[
\left\Vert u_{1}-u_{2}\right\Vert _{Y^{s,p}}\left( \left\Vert
u_{1}\right\Vert _{X_{\infty }}+\left\Vert u_{2}\right\Vert _{X_{\infty
}}\right) ^{\gamma }+ 
\]%
\[
\left. \left( \left\Vert u_{1}\right\Vert _{X_{\infty }}+\left\Vert
u_{2}\right\Vert _{X_{\infty }}\right) ^{\gamma -1}\left\Vert
u_{1}+u_{2}\right\Vert _{X_{p}}\left\Vert u_{1}-u_{2}\right\Vert
_{X_{p}}\right\} \leq 
\]%
\[
C\left( \left\Vert u_{1}\right\Vert _{C\left( T,A\right) }+\left\Vert
u_{2}\right\Vert _{C\left( T,A\right) }\right) ^{\gamma }\left\Vert
u_{1}-u_{2}\right\Vert _{C\left( T,A\right) }. 
\]

On the oher hand, by $(2.17)$ in Theorem 2.2, Corollary 3.3 and $\left(
3.7\right) $, we get

\[
\left( \left\Vert A^{\alpha }\ast \left[ G\left( u_{1}\right) -G\left(
u_{2}\right) \right] \right\Vert _{Y^{s,p}}+\left\Vert A^{\alpha }\ast \left[
G\left( u_{1}\right) -G\left( u_{2}\right) \right] _{t}\right\Vert
_{Y^{s,p}}\right) \leq 
\]

\begin{equation}
C\dint\limits_{0}^{t}\left\Vert \Delta g\ast \left[ f\left( u_{1}\right)
\left( \tau \right) -f\left( u_{2}\right) \left( \tau \right) \right]
\right\Vert _{Y^{s,p}}d\tau \leq  \tag{3.15}
\end{equation}%
\[
C\dint\limits_{0}^{t}\left\{ \left\Vert u_{1}-u_{2}\right\Vert _{X_{\infty
}}\left( \left\Vert u_{1}\right\Vert _{Y^{s,p}}+\left\Vert u_{2}\right\Vert
_{Y^{s,p}}\right) \left( \left\Vert u_{1}\right\Vert _{X_{\infty
}}+\left\Vert u_{2}\right\Vert _{X_{\infty }}\right) ^{\gamma -1}+\right. 
\]%
\[
\left. \left\Vert u_{1}-u_{2}\right\Vert _{Y^{s,p}}\left( \left\Vert
u_{1}\right\Vert _{X_{\infty }}+\left\Vert u_{2}\right\Vert _{X_{\infty
}}\right) ^{\gamma }\right\} d\tau \leq 
\]%
\[
C\left( \left\Vert u_{1}\right\Vert _{C\left( T,A\right) }+\left\Vert
u_{2}\right\Vert _{C\left( T,A\right) }\right) ^{\gamma }\left\Vert
u_{1}-u_{2}\right\Vert _{C\left( T,A\right) }. 
\]

Combining $\left( 3.14\right) $ with $\left( 3.15\right) $ yields 
\begin{equation}
\left\Vert G\left( u_{1}\right) -G\left( u_{2}\right) \right\Vert _{C\left(
T,A\right) }\leq  \tag{3.16}
\end{equation}

\[
C\left( \left\Vert u_{1}\right\Vert _{C\left( T,A\right) }+\left\Vert
u_{2}\right\Vert _{C\left( T,A\right) }\right) ^{\gamma }\left\Vert
u_{1}-u_{2}\right\Vert _{C\left( T,A\right) }. 
\]

\bigskip Taking $\delta $ is enough small, from $\left( 3.16\right) $ we
obtain that $G$ is strictly controctive in $C\left( T,A\right) $. Using the
contaction mapping principle we get that $G\left( u\right) $ has a unique
fixed point $u\left( x,t\right) \in C\left( T,A\right) $ and $u\left(
x,t\right) $ is the solution $\left( 1.1\right) -\left( 1.2\right) .$

\bigskip Let us show that this solution is a unique in $C^{2,s}\left(
A,E\right) $. Let $u_{1}$, $u_{2}\in C^{2,s}\left( A,E\right) $ are two
solution of $(1.1)-(1.2)$. Then for $u=u_{1}-u_{2}$, we have 
\begin{equation}
u_{tt}-a\ast \Delta u+A\ast u=\Delta g\ast \left[ f\left( u_{1}\right)
-f\left( u_{2}\right) \right] \text{ }  \tag{3.17}
\end{equation}%
Hence, by Minkowski's inequality for integrals and by Theorem 2.2 from $%
\left( 3.17\right) $ we obtain

\begin{equation}
\left\Vert u_{1}-u_{2}\right\Vert _{Y^{s,p}}\leq C_{2}\left( T\right) \text{ 
}\dint\limits_{0}^{t}\left\Vert u_{1}-u_{2}\right\Vert _{Y^{s,p}}d\tau . 
\tag{3.18}
\end{equation}%
From $(3.18)$ and Gronwall's inequality, we have $\left\Vert
u_{1}-u_{2}\right\Vert _{Y^{s,p}}=0$, i.e. problem $(1.1)-(1.2)$ has a
unique solution in $C^{2,s}\left( A,E\right) $.

Consider the problem $\left( 1.1\right) -\left( 1.2\right) $, when $\varphi
\in \mathbb{E}_{0p}$ and $\psi \in \mathbb{E}_{1p}$. Let 
\[
C^{\left( i\right) }\left( Y^{s,2}\right) =C^{\left( i\right) }([0,\infty
);Y^{s,2}\left( A,E\right) )\text{, }i=0,1,2. 
\]

By reasoning as in Theorem 3.1 and $\left[ \text{13, Theorem 1.1}\right] $
we have:

\textbf{Condition 3.2. }Assume: (1) the Condition 2.1 holds; (2) $0<\alpha
<1-\frac{1}{2p}$, $\varphi \in $ $\mathbb{E}_{0p}$, $\psi \in $ $\mathbb{E}%
_{1p}$ and $s>\frac{n}{p}$; (3) $f\in C^{\left[ s\right] +1}\left( \mathbb{R}%
;E\right) $ with $f(0)=0$; (4) the kernel $g$ is a bounded operator function
in $E$, whose Fourier transform satisfies 
\begin{equation}
0\leq \left\Vert \hat{g}\left( \xi \right) \right\Vert _{B\left( E\right)
}\leq C_{g}\left( 1+\left\vert \xi \right\vert ^{2}\right) ^{-\frac{r}{2}}%
\text{ for all }\xi \in \mathbb{R}^{n}\text{ and }r\geq 2.  \tag{3.19}
\end{equation}

\ \textbf{Theorem 3.2. }Let the Condition 3.2 holds. Then there exists a
constant $\delta >0$, such that for any $\varphi \in $ $\mathbb{E}_{0p}$, $%
\psi \in $ $\mathbb{E}_{1p}$ satisfying 
\begin{equation}
\left\Vert \varphi \right\Vert _{\mathbb{E}_{0p}}+\left\Vert \psi
\right\Vert _{\mathbb{E}_{1p}}\leq \delta ,  \tag{3.20}
\end{equation}%
problem $\left( 1.1\right) -\left( 1.2\right) $ has a unique local strange
solution $u\in C^{\left( 2\right) }\left( Y^{s,p}\right) $. Moreover,

\begin{equation}
\sup_{t\in \left[ 0,T\right] }\left( \left\Vert u\left( .,t\right)
\right\Vert _{\hat{Y}^{s,p}\left( A^{\alpha },E\right) }+\left\Vert
u_{t}\left( .,t\right) \right\Vert _{\hat{Y}^{s,p}\left( A^{\alpha
};E\right) }\right) \leq C\delta ,  \tag{3.21}
\end{equation}%
where the constant $C$ only depends on $f$ \ and initial data.

\textbf{Proof. }Consider a metric space defined by 
\[
W_{0}^{s,p}=\left\{ u\in C^{\left( 2\right) }\left( Y^{s,p}\right) \text{, }%
\left\Vert u\right\Vert _{Y^{s,p}}\leq 3C_{0}\delta \right\} , 
\]%
equipped with the norm%
\[
\left\Vert u\right\Vert _{W_{0}^{s,p}}=\sup\limits_{t\in \left[ 0,T\right]
}\left( \left\Vert u\right\Vert _{\hat{Y}^{s,p}\left( A^{\alpha };E\right)
}+\left\Vert u_{t}\right\Vert _{\hat{Y}^{s,p}\left( A^{\alpha };E\right)
}\right) , 
\]%
where $\delta >0$ satisfies $\left( 3.20\right) $ and $C_{0}$ \ is a
constant in Theorem 2.1. It is easy to prove that $W_{0}^{s,p}$ is a
complete metric space. From Sobolev imbedding theorem we know that $%
\left\Vert u\right\Vert _{\infty }\leq 1$ if we take that $\delta $ is
enough small. By Theorem 2.2 and Corollary 3.1, $\Delta g\ast f\left(
u\right) \in Y^{s,p}$. Thus the problem $\left( 1.1\right) -\left(
1.2\right) $ has a unique solution which can be written as $\left(
3.9\right) $. We should prove that the operator $G\left( u\right) $\ defined
by $\left( 3.9\right) $ is strictly contractive if $\delta $ is suitable
small. In fact, by $(2.17)$ in Theorem 2.2 and Lemma 3.1 we get 
\[
\left\Vert A^{\alpha }\ast G\left( u\right) \right\Vert
_{Y^{s,p}}+\left\Vert A^{\alpha }\ast G_{t}\left( u\right) \right\Vert
_{Y^{s,p}}\leq C_{0}\left[ \left\Vert \varphi \right\Vert _{\mathbb{E}%
_{0p}}\right. +\left\Vert \psi \right\Vert _{\mathbb{E}_{1p}}+ 
\]

\[
\left. \dint\limits_{0}^{t}\left\Vert K\left( u\right) \left( .,\tau \right)
\right\Vert _{Y^{s,p}}d\tau \right] \leq C_{0}\delta
+C_{0}\dint\limits_{0}^{t}\left\Vert K\left( u\right) \left( .,\tau \right)
\right\Vert _{Y^{s,p}}d\tau \leq 
\]%
\begin{equation}
C_{0}\delta +C\dint\limits_{0}^{t}\left\Vert u\left( \tau \right)
\right\Vert _{Y^{s,p}}d\tau \leq C_{0}\delta +C\left\Vert u\right\Vert
_{Y^{s,p}},  \tag{3.22}
\end{equation}%
where 
\[
K\left( u\right) \left( .,\tau \right) =S\left( x,t-\tau \right) \Delta
g\ast f\left( u\right) \left( x,\tau \right) . 
\]%
Therefore, from $(3.22)$ we have 
\begin{equation}
\left\Vert G\left( u\right) \right\Vert _{Y^{s,p}}\leq 2C_{0}\delta
+C\left\Vert u\right\Vert _{Y^{s,p}}\text{.}  \tag{3.23}
\end{equation}%
Taking that $\delta $ is enough small such that $C\left( 3C_{0}\delta
\right) ^{\alpha }$ $<1/3$, from $\left( 3.23\right) $ and from Theorems
2.1, 2.2 we get that $G$ maps $W_{0}^{s,p}$ into $W_{0}^{s,p}$. Then, by
reasoning as in Theorem 3.1 we obtain that $G$ : $W_{0}^{s,p}\rightarrow
W_{0}^{s,p}$ is strictly contractive. Using the contraction mapping
principle, we know that $G(u)$ has a unique fixed point $u\in $ $C^{\left(
2\right) }\left( Y^{s,2}\right) $ and $u(x,t)$ is the solution of the
problem $(1.1)-(1.2)$. Moreover, by virtue of Theorem 2.1 from $\left(
3.20\right) $ we obtain $\left( 3.21\right) .$

We claim that the solution of $(1.1)-(1.2)$ is also unique in $C^{\left(
1\right) }\left( Y^{s,2}\right) $. In fact, let $u_{1}$ and $u_{2}$ be two
solutions of the problem $(1.1)-(1.2)$ and $u_{1}$, $u_{2}\in C^{\left(
2\right) }\left( Y^{s,2}\right) $. Using the contraction mapping principle,
we know that $G(u)$ has a unique fixed point $u\in C^{\left( 2\right)
}\left( Y^{s,2}\right) $. Using the contraction mapping principle, we know
that $G(u)$ has a unique fixed point $u\in $ $C^{\left( 2\right) }\left(
Y^{s,2}\right) $. Let $u=u_{1}-u_{2}$, then%
\[
u_{tt}-a\ast \Delta u+A\ast u=\Delta \left[ g\ast \left( f\left(
u_{1}\right) -f\left( u_{2}\right) \right) \right] . 
\]

\bigskip This fact is derived in a similar way as in Theorem 3.1, by using
Theorems 2.1, 2.2 and Gronwall's inequality.

Let 
\[
C^{\left( 2,s\right) }\left( Y^{s,p}\right) =C^{\left( 2\right) }\left( 
\left[ 0,T\right] ;Y^{s,p}\left( A;E\right) \right) . 
\]

\textbf{Theorem 3.3. }Let the Condition 3.2 hold. Then there is some $T>0$
such that the problem $(1.1)-(1.2)$ for initial data $\varphi \in $ $\mathbb{%
E}_{0p}$ and $\psi \in $ $\mathbb{E}_{1p}$ is well posed with solution in $%
C^{1}\left( \left[ 0,T\right] ;Y^{s,p}\left( A,E\right) \right) .$

\textbf{Proof. }Consider the convolution operator 
\[
u\rightarrow \Delta \left[ g\ast f\left( u\right) \right] . 
\]%
In view of assumptions and Fourier multipler results in $X_{p}$ spaces (see
e.g. $\left[ \text{12, Theorem 4.3}\right] $)\ we have \ 
\[
\left\Vert \Delta g\ast \upsilon \right\Vert _{Y^{s,p}}\lesssim \left\Vert 
\mathbb{F}^{-1}\left( 1+\xi \right) ^{\frac{s}{2}}\left\vert \xi \right\vert
^{2}\hat{g}\left( \xi \right) \hat{\upsilon}\left( \xi \right) \right\Vert
\lesssim \left\Vert \upsilon \right\Vert _{Y^{s,p}}, 
\]%
i.e. $\Delta g\ast \upsilon $ is a bounded linear operator on $Y^{s,p}$.
Then by Corollary 3.1, $K\left( u\right) $ is locally Lipschitz on $Y^{s,p}$%
. Then by reasoning as in Theorem 3.2 and $\left[ \text{13, Theorem 1.1}%
\right] $ we obtain that $G$: $W_{0}^{s,p}\rightarrow W_{0}^{s,p}$ is
strictly contractive. Using the contraction mapping principle, we get that
the operator $G(u)$ defined by $\left( 3.5\right) $ has a unique fixed point 
$u(x,t)\in $ $C^{\left( 2\right) }\left( Y^{s,p}\right) $ \ and $u(x,t)$ is
the solution of the problem $(1.1)-(1.2)$. Moreover, we show that the
solution $u(x,t)$ of $(1.1)-(1.2)$ is also unique in $C^{\left( 2\right)
}\left( Y^{s,p}\right) $. In fact, let $u_{1}$ and $u_{2}$ be two solutions
of the problem $(1.1)-(1.2)$ and $u_{1}$, $u_{2}\in C^{\left( 2\right)
}\left( Y^{s,p}\right) $. Let $u=u_{1}-u_{2}$, then%
\[
u_{tt}-a\ast \Delta u+A\ast u=\Delta \left[ g\ast \left( f\left(
u_{1}\right) -f\left( u_{2}\right) \right) \right] . 
\]%
This fact is derived in a similar way as in Theorem 3.2, by using Theorems
2.1, 2.2 and Gronwall's inequality.

\textbf{Theorem 3.4. }Let the Condition 3.2 holds for $r>2+\frac{n}{p}$.
Then there is some $T>0$ such that problem $(1.1)-(1.2)$ is well posed for $%
\varphi \in $ $\mathbb{E}_{0p}$ and $\psi \in $ $\mathbb{E}_{1p}$ with
solution in $C^{\left( 2\right) }\left( Y^{s,p}\right) .$

\textbf{Proof. } All we need here, is to show that $K\ast f(u)$ is Lipschitz
on $Y^{s,p}$. Indeed, by reasoning as in Theorem 3.3 we have 
\[
\left\Vert \Delta g\ast \upsilon \right\Vert _{Y^{s+r-p}}\lesssim \left\Vert
\left( 1+\left\vert \xi \right\vert ^{2}\right) ^{\frac{s+r-2}{2}}\left\vert
\xi \right\vert ^{2}\hat{g}\left( \xi \right) \hat{\upsilon}\left( \xi
\right) \right\Vert \lesssim \left\Vert \upsilon \right\Vert _{Y^{s,p}}, 
\]

\bigskip Then $\Delta g\ast \upsilon $ is a bounded linear map from $Y^{s,p}$
into $Y^{s+r-2,p}$. Since $s\geq 0$ and $r$ $>2+\frac{n}{p}$ \ we get 
\[
s+r-2>\frac{n}{p}. 
\]
The embedding theorem for $E-$valued Sobolev spaces (see e.g, $\left[ 31%
\right] $) implies that $\Delta g\ast \upsilon $\ is a bounded linear map
from $Y^{s,p}\left( A;E\right) $ into $Y^{s,p}\left( A;E\right) $. Lemma 3.2
implies the Lipschitz condition on $Y^{s,p}$. Then, by reasoning as in
Theorem 3.3 we obtain the assertion.

The solution in theorems 3.2-3.4 can be extended to a maximal interval $%
[0,T_{\max })$, where finite $T_{\max }$ is characterized by the blow-up
condition 
\[
\limsup\limits_{T\rightarrow T_{\max }}\left\Vert u\right\Vert _{\hat{Y}%
^{s,p}\left( A^{\alpha };E\right) }=\infty . 
\]

\textbf{Lemma 3.8.} \textbf{\ }Let the Condition 3.2 hold and $u$ is a
solution of $(1.1)-(1.2)$. Then there is a global solution if for any $%
T<\infty $ we have%
\begin{equation}
\sup\limits_{t\in \left[ 0,T\right] }\left( \left\Vert u\right\Vert _{\hat{Y}%
^{s,p}\left( A^{\alpha };E\right) }+\left\Vert u_{t}\right\Vert _{\hat{Y}%
^{s,p}\left( A^{\alpha };E\right) }\right) <\infty .  \tag{3.24}
\end{equation}

\textbf{Proof. }Indeed, by reasoning as in the second part of the proof of
Theorem 3.1, by using a continuation of local solution of $(1.1)-(1.2)$ and
assuming contrary that, $\left( 3.24\right) $ holds and $T_{0}<\infty $\ \
we obtain contradiction, i.e. we get $T_{0}=T_{\max }=\infty .$

\begin{center}
\textbf{4. Conservation of energy and global existence. }
\end{center}

In this section, we prove the existence and the uniqueness of the global
strong solution for the problem $(1.1)-(1.2).$ \ For this purpose, we are
going to make a priori estimates of the local strong solution of $%
(1.1)-(1.2).$

\textbf{Condition 4.1. }Suppose the Condition 3.2 is satisfied.\textbf{\ }%
Assume $a\in L^{2}\left( \mathbb{R}^{n}\right) $ and the kernel $g$ is a
bounded operator function in $E$, whose Fourier transform satisfies%
\[
0<\left\Vert \hat{g}\left( \xi \right) \right\Vert _{B\left( E\right)
}\lesssim \left( 1+\left\vert \xi \right\vert ^{2}\right) ^{-\frac{r}{2}}%
\text{ for all }\xi \in \mathbb{R}^{n}\text{ and }r\leq 2\left( s+1\right) . 
\]%
Moreover, let $\hat{g}\left( \xi \right) $ have fractional powers for all $%
\xi \in \mathbb{R}^{n}$. Let $\mathbb{F}^{-1}$ denote the inverse Fourier
transform. Assume that the operator $\hat{g}\left( \xi \right) $ has a
fractinal pover $\hat{g}^{^{\frac{1}{2}}}\left( \xi \right) $ for all $\xi
\in \mathbb{R}^{n}$. We consider the Fourier multipler operator $B=B_{g}$
defined by 
\[
u\in D\left( B\right) =Y^{s,p}\text{, }Bu=\mathbb{F}^{-1}\left[ \left\vert
\xi \right\vert ^{-1}\hat{g}^{^{-\frac{1}{2}}}\left( \xi \right) \hat{u}%
\left( \xi \right) \right] , 
\]%
Then it is clear to see that 
\begin{equation}
B^{-2}u=-\Delta g\ast u\text{, }B^{-1}u=\mathbb{F}^{-1}\left[ \left\vert \xi
\right\vert \hat{g}^{\frac{1}{2}}\left( \xi \right) \hat{u}\left( \xi
\right) \right] .  \tag{4.1}
\end{equation}%
Let 
\[
C^{\left( 1\right) }\left( L^{p}\right) =C^{\left( 1\right) }\left( \left[
0,\right. \left. T\right) ;L^{p}\left( \mathbb{R}^{n};E\right) \right) \text{%
, }C^{\left( 2,s\right) }\left( A,E\right) =C^{\left( 2\right) }\left( \left[
0,T\right] ;Y^{s,p}\left( A;E\right) \right) , 
\]%
where $Y^{s,p}\left( A;E\right) $ was defined in Section 2.

First, we show the following\bigskip

\textbf{Lemma 4.1. }Let the Condition 4.1 holds and $0<\alpha <1-\frac{1}{2p}
$. Assume there exist a solution $u\in C^{\left( 2,s\right) }\left(
A,E\right) $ of $(1.1)-(1.2)$. Then%
\[
\hat{A}^{\alpha }Bu\text{, }\hat{A}^{\alpha }Bu_{t}\in C^{\left( 1\right)
}\left( L^{p}\right) . 
\]

\bigskip \textbf{Proof.} By Lemma 2.1, problem $\left( 1.1\right) -\left(
1.2\right) $ is equivalent to the followng integra equation,%
\begin{equation}
u\left( x,t\right) =C_{1}\left( t\right) \varphi +S_{1}\left( t\right) \psi
+Qg,  \tag{4.2}
\end{equation}%
where $C_{1}\left( t\right) $, $S_{1}\left( t\right) $ are operator
functions defined by $\left( 2.5\right) $ and $\left( 2.6\right) $, where $g$
replaced by $g\ast f\left( u\right) $ and 
\begin{equation}
Qg=\dint\limits_{0}^{t}\mathbb{F}^{-1}\left[ S\left( \xi ,t-\tau \right)
\left\vert \xi \right\vert ^{2}\hat{g}\left( \xi \right) \hat{f}\left(
u\right) \left( \xi \right) \right] d\tau .  \tag{4.3}
\end{equation}%
From $\left( 4.2\right) $ we get that 
\[
u_{t}\left( x,t\right) =\frac{d}{dt}C_{1}\left( t\right) \varphi +\frac{d}{dt%
}S_{1}\left( t\right) \psi + 
\]%
\begin{equation}
\dint\limits_{0}^{t}\mathbb{F}^{-1}\left[ C\left( \xi ,t-\tau \right)
\left\vert \xi \right\vert ^{2}\hat{g}\left( \xi \right) \hat{f}\left(
G\left( u\right) \left( \xi \right) \right) \right] d\tau .  \tag{4.4}
\end{equation}%
Since $C_{1}\left( t\right) ,$ $S_{1}\left( t\right) $ and $\frac{d}{dt}%
S\left( \xi ,t\right) $\ are uniformly bounded operators in $E$ for fixet $t$%
, by $\left( 4.1\right) $, $\left( 4.2\right) $, $\left( 4.4\right) $ and
Fourier multipler results in $X_{p}$ spaces (see e.g. $\left[ \text{12,
Theorem 4.3}\right] $) we have%
\begin{equation}
\left\Vert \hat{A}^{\alpha }BC_{1}\left( t\right) \varphi \right\Vert
_{L^{p}}=\left\Vert \mathbb{F}^{-1}\left[ \left\vert \xi \right\vert ^{-1}%
\hat{g}^{-\frac{-1}{2}}\left( \xi \right) \hat{A}^{\alpha }C\left( \xi
,t\right) \hat{\varphi}\right] \right\Vert _{L^{p}}\lesssim  \tag{4.5}
\end{equation}%
\[
\left\Vert \varphi \right\Vert _{\mathbb{E}_{0p}}<\infty , 
\]%
\[
\left\Vert \hat{A}^{\alpha }BS_{1}\left( t\right) \varphi \right\Vert
_{L^{p}}=\left\Vert \mathbb{F}^{-1}\left[ \left\vert \xi \right\vert ^{-1}%
\hat{g}^{-\frac{-1}{2}}\left( \xi \right) \hat{A}^{\alpha }S\left( \xi
,t\right) \hat{\psi}\right] \right\Vert _{L^{p}}\lesssim 
\]%
\[
\left\Vert \psi \right\Vert _{\mathbb{E}_{1p}}<\infty . 
\]%
By differentiating $\left( 2.3\right) $, in a similar way we have 
\[
\left\Vert \hat{A}^{\alpha }B\frac{d}{dt}C_{1}\left( t\right) \varphi
\right\Vert _{L^{p}}=\left\Vert F^{-1}\left[ \left\vert \xi \right\vert ^{-1}%
\hat{g}^{-\frac{-1}{2}}\left( \xi \right) \hat{A}^{\alpha }\frac{d}{dt}%
C\left( \xi ,t\right) \hat{\varphi}\right] \right\Vert _{L^{p}} 
\]%
\begin{equation}
\lesssim \left\Vert \varphi \right\Vert _{\mathbb{E}_{0p}}<\infty , 
\tag{4.6}
\end{equation}%
\[
\left\Vert \hat{A}^{\alpha }B\frac{d}{dt}S_{1}\left( t\right) \varphi
\right\Vert _{L^{p}}=\left\Vert F^{-1}\left[ \left\vert \xi \right\vert ^{-1}%
\hat{g}^{-\frac{-1}{2}}\left( \xi \right) \hat{A}^{\alpha }\frac{d}{dt}%
S\left( \xi ,t\right) \hat{\psi}\right] \right\Vert _{L^{p}}\lesssim 
\]%
\[
\left\Vert \psi \right\Vert _{\mathbb{E}_{1p}}<\infty . 
\]%
For fixed $t$, we have $f(u)\in Y^{s,p}$. Moreover, by assumption on $\hat{A}%
\left( \xi \right) $ we have the uniformly estimate 
\[
\left\Vert \hat{A}^{\alpha }\left( \xi \right) \eta ^{-1}\left( \xi \right)
\right\Vert _{B\left( E\right) }\leq C_{A}. 
\]%
\ Then by hypothesis on $\hat{g}\left( \xi \right) $, due to $s+r\geq 1$\
from $\left( 4.1\right) $ and $\left( 4.3\right) $ and Fourier multipler
results in $X_{p}$\ we get 
\[
\left\Vert \hat{A}^{\alpha }BQg\right\Vert _{L^{p}}\leq \left\Vert \mathbb{F}%
^{-1}\left[ \left\vert \xi \right\vert \hat{g}^{\frac{1}{2}}\left( \xi
\right) \hat{A}^{\alpha }\left( \xi \right) \dint\limits_{0}^{t}S\left( \xi
,t-\tau \right) \hat{f}\left( u\right) \left( \xi \right) d\tau \right]
\right\Vert _{L^{p}}\lesssim 
\]%
\begin{equation}
C_{A}\left\Vert f\left( u\right) \right\Vert _{Y^{s,p}}<\infty .  \tag{4.7}
\end{equation}%
Then from $\left( 4.2\right) $ and $\left( 4.4\right) -\left( 4.7\right) $
we obtain the assertion.

\textbf{Lemma 4.2.} Assume the Condition 4.1 holds with $a=0$. Moreover, let%
\[
\left\Vert \left( \hat{g}\left( \xi \right) \right) ^{-\frac{-1}{2}%
}\right\Vert _{B\left( E\right) }=O\left( 1+\left\vert \xi \right\vert
^{2}\right) ^{\frac{s+1}{2}}. 
\]%
Suppose the solution of $(1.1)-(1.2)$ exists in $C^{\left( 2,s\right)
}\left( A,E\right) $. If $B\psi \in L^{p}$ then $Bu_{t}\in C^{\left(
1\right) }\left( L^{p}\right) $. Moreover, if $B\varphi \in L^{p}$, then $%
Bu\in C^{\left( 1\right) }\left( L^{p}\right) .$

\textbf{Proof. }Integrating the equation $\left( 1.1\right) $ for $a=0$
twice and calculating the resulting double integral as an iterated integral,
we have 
\begin{equation}
u\left( x,t\right) =\varphi \left( x\right) +t\psi \left( x\right) - 
\tag{4.8}
\end{equation}%
\[
\dint\limits_{0}^{t}\left( t-\tau \right) \left( A\ast u\right) \left(
x,\tau \right) d\tau +\dint\limits_{0}^{t}\left( t-\tau \right) \Delta
\left( g\ast f\left( u\right) \right) \left( x,\tau \right) d\tau ,
\]

\begin{equation}
u_{t}\left( x,t\right) =\psi \left( x\right) -\dint\limits_{0}^{t}\left(
A\ast u\right) \left( x,\tau \right) d\tau +\dint\limits_{0}^{t}\Delta
\left( g\ast f\left( u\right) \right) \left( x,\tau \right) d\tau . 
\tag{4.9}
\end{equation}%
From $\left( 4.1\right) $ and $\left( 4.9\right) $ for fixed $t$ and $\tau $
we get $f\left( u\right) \in Y^{s,p}$ for all $t$. Also 
\begin{equation}
\left\Vert B\Delta \left( g\ast f\left( u\right) \right) \left( x,\tau
\right) \right\Vert _{L^{p}}\lesssim \left\Vert \mathbb{F}^{-1}\left[
\left\vert \xi \right\vert \left( \hat{g}^{\frac{-1}{2}}\left( \xi \right)
\right) \hat{f}\left( u\right) \left( \xi \right) \right] \right\Vert
_{L^{p}}.  \tag{4.10}
\end{equation}%
Then from $\left( 4.8\right) -\left( 4.10\right) $ we obtain%
\[
\left\Vert Bu_{t}\left( x,t\right) \right\Vert _{L^{2}}=\left\Vert B\psi
\left( x\right) \right\Vert _{L^{2}}- 
\]%
\[
\dint\limits_{0}^{t}\left\Vert B\left( A\ast u\right) \left( x,\tau \right)
\right\Vert _{L^{p}}d\tau -\dint\limits_{0}^{t}\left\Vert B\Delta \left(
g\ast f\left( u\right) \right) \left( x,\tau \right) \right\Vert
_{L^{p}}d\tau . 
\]%
By assumption on $A$, $g$ and by $\left( 4.1\right) $ for fixed $\tau $ we
have $Bu_{t}\in C^{\left( 1\right) }\left( L^{p}\right) .$%
\[
\left\Vert B\left( A\ast u\right) \left( x,\tau \right) \right\Vert
_{L^{p}}\lesssim \left\Vert \mathbb{F}^{-1}\left[ \left\vert \xi \right\vert
^{-1}\hat{A}\left( \xi \right) \left( \hat{g}^{-\frac{-1}{2}}\left( \xi
\right) \right) \hat{u}\left( \xi ,\tau \right) \right] \right\Vert _{L^{p}} 
\]%
\[
\lesssim \left\Vert u\left( .,\tau \right) \right\Vert _{Y^{s,p}\left(
A\right) }. 
\]%
Moreover, by Lemma 3.3 we have $Bu_{t}\in C^{\left( 1\right) }\left(
L^{p}\right) $. The second statement follows similarly from $\left(
4.8\right) .$

From Lemma 4.2 we obtain the following result.

\textbf{Result 4.1.} Assume the Condition 4.1 are satisfied with $a=0$ and%
\[
\left\Vert \hat{g}\left( \xi \right) \right\Vert _{B\left( E\right)
}=O\left( 1+\left\vert \xi \right\vert ^{2}\right) ^{-\frac{r}{2}}. 
\]%
Suppose the solution of $(1.1)-(1.2)$ exists in $C^{\left( 2,s\right)
}\left( A,H\right) $ for some $s\geq 0$. If $B\psi \in L^{2}$ then $%
Bu_{t}\in C^{\left( 1\right) }\left( L^{2}\right) $. Moreover, if $B\varphi
\in L^{2}$, then $Bu\in C^{\left( 1\right) }\left( L^{2}\right) .$

\textbf{Lemma 4.3. }Assume the Condition 3.2 holds and $s+r\geq 1$. Let $%
u\in C^{\left( 2,s\right) }\left( A,H\right) $ be a slutlion of $\left(
1.1\right) -\left( 1.2\right) $ for any $t\in \left[ 0,\right. \left.
T\right) $. Let $B\psi \in L^{2}$ and $\left( f\left( u\right) ,u\right) \in
L^{2}$. Then the energy 
\begin{equation}
E\left( t\right) =\left\Vert Bu_{t}\right\Vert _{L^{2}}^{2}+\left( B\left[
A\ast u-a\ast \Delta u\right] ,Bu\right) _{L^{2}}+\left( f\left( u\right)
,u\right) _{L^{2}}  \tag{4.11}
\end{equation}%
is constant.

\textbf{Proof. }By Theorem 4.1, $\ A^{\alpha }Bu$, $A^{\alpha }Bu_{t}\in
L^{2}$ for $0<\alpha <\frac{3}{4}$. By assumptions $\left( f\left( u\right)
,u\right) \in L^{2}$ and $A\ast u\in L^{2}$. By use of $\left( 1.1\right) $
and Parseval's identity, it follows from straightforward calculation that%
\[
\frac{d}{dt}E\left( t\right) =2\left( Bu_{tt},Bu_{t}\right) +2\left( Ba\ast
\Delta u,Bu_{t}\right) + 
\]

\[
2\left[ B\left( A\ast u\right) ,Bu_{t}\left( t\right) \right] +2\left(
f\left( u\right) ,u_{t}\right) =2B^{2}\left( u_{tt},u_{t}\right) + 
\]%
\[
2B^{2}\left( a\ast \Delta u,u_{t}\right) +2B^{2}\left( A\ast u,u_{t}\right)
-2B^{2}\left( \Delta \left[ g\ast f\left( u\right) \right] ,u_{t}\right) = 
\]%
\[
2B^{2}\left( u_{tt}+a\ast \Delta u+A\ast u-\Delta \left[ g\ast f\left(
u\right) \right] ,u_{t}\right) = 
\]%
\[
B^{2}\frac{d}{dt}\left[ \left( u_{tt}-a\ast \Delta u+A\ast u-\Delta g\ast
f\left( u\right) ,u\right) \right] =0, 
\]%
where $\left( u,\upsilon \right) $ denotes the inner product in $L^{2}\left( 
\mathbb{R}^{n}\right) $. Hence, we obtain the assertion.

By using the above lemmas we obtain the following results

\textbf{Theorem 4.1. }Let the Condition 4.1 holds for $r>2+\frac{n}{2}$.
Moreover, let $B\psi \in L^{2}$, $\left( f\left( u\right) ,u\right) \in
L^{2}\left( \mathbb{R}^{n};H\right) $ and there is some $k>0$ so that $%
\left( f\left( u\right) ,u\right) \geq -k\left\Vert u\left( .,t\right)
\right\Vert ^{2}$ for all $t\in \left[ 0,T\right] $. Then there is some $T>0$
such that problem $(1.1)-(1.2)$ has a global solution $u\in C^{\left(
2,s\right) }\left( A,H\right) .$ \ \ \ 

\textbf{Proof. }Since $r>2+\frac{n}{2}$, by Theorem 3.4 we get local
existence $u\in C^{\left( 2,s\right) }\left( A,E\right) $ for some $T>0$.
Assume that $u$ exists on $[0,T)$. By Lemma 4.3,\ we obtain%
\begin{equation}
\left\Vert Bu_{t}\right\Vert ^{2}+\left\Vert \hat{a}\right\Vert
_{L^{2}}\left\Vert \mathbb{F}^{-1}\hat{g}\ast \hat{u}\right\Vert ^{2}+\left(
B\left( A\ast u\right) ,Bu\right) \leq   \tag{4.12}
\end{equation}%
\[
E\left( 0\right) +2k\left\Vert u\left( .,t\right) \right\Vert ^{2}.
\]%
Let $Y^{s,2}$ denotes by $W^{s}$. By condition on $\hat{g}\left( \xi \right) 
$, we have 
\begin{equation}
\left\Vert Bu_{t}\right\Vert _{L^{2}\left( A\right) }^{2}=\dint\limits_{%
\mathbb{R}^{n}}\left\vert \xi \right\vert ^{-2}\left\Vert \hat{g}^{-1}\left(
\xi \right) \right\Vert _{B\left( H\right) }^{2}\left\Vert A\hat{u}%
_{t}\left( \xi ,t\right) \right\Vert _{H}^{2}\geq   \tag{4.13}
\end{equation}%
\[
C_{g}^{-1}\dint\limits_{\mathbb{R}^{n}}\left( 1+\left\vert \xi \right\vert
^{2}\right) ^{\frac{r}{2}-1}\left\Vert A\hat{u}_{t}\left( \xi ,t\right)
\right\Vert _{H}^{2}\approx C_{g}^{-1}\left\Vert Au_{t}\left( t\right)
\right\Vert _{W^{\frac{r}{2},-1}}^{2},
\]%
where $C_{g}$ is the positive constant that appears in $\left( 3.19\right) $%
. By properties of norms in Hilbert spaces and by Cauchy-Schwarz inequality,
from $\left( 4.12\right) $ and $\left( 4.13\right) $ we get 
\[
\frac{d}{dt}\left\Vert u\left( t\right) \right\Vert _{W^{\frac{r}{2}%
-1}\left( A\right) }^{2}\leq 2\left\Vert u_{t}\left( t\right) \right\Vert
_{W^{\frac{r}{2}-1}\left( A\right) }\left\Vert u\left( t\right) \right\Vert
_{W^{\frac{r}{2}-1}\left( A\right) }\leq 
\]%
\[
\left\Vert u_{t}\left( t\right) \right\Vert _{W^{\frac{r}{2}-1}\left(
A\right) }^{2}+\left\Vert u\left( t\right) \right\Vert _{W^{\frac{r}{2}%
-1}\left( A\right) }^{2}\leq C\left\Vert Bu_{t}\left( t\right) \right\Vert
_{W^{\frac{r}{2}-1}\left( A\right) }^{2}+
\]%
\[
\left\Vert u\left( t\right) \right\Vert _{W^{\frac{r}{2}-1}\left( A\right)
}^{2}\leq CE\left( 0\right) +\left( 2Ck+1\right) \left\Vert u\left( t\right)
\right\Vert _{W^{\frac{r}{2}-1}\left( A\right) }^{2}\text{.}
\]

Gronwall's lemma implies that $\left\Vert u\left( t\right) \right\Vert _{W^{%
\frac{r}{2}-1}\left( A\right) }$ is bounded in $[0,T)$. But, since $\frac{r}{%
2}-1>\frac{n}{4}$, we conclude that $\left\Vert u(t)\right\Vert _{L^{\infty
}\left( A\right) }$ also is bounded in $[0,T)$. By Lemma 3.8 this implies a
global solution.

\begin{center}
\textbf{5. Blow up in finite time}
\end{center}

We will use the following lemma to prove blow up in fininite time.

\textbf{Lemma 5.1} $[$16$]$ Suppose $H(t)$, $t\geq 0$ is a positive, twice
differentiable function satisfying $H^{\left( 2\right) }H-\left( 1+\nu
\right) \left( H^{\left( 1\right) }\right) ^{2}\geq 0$, where $\nu >$ $0$.
If $H(0)$ $>0$ and $H^{\left( 1\right) }(0)>0$, then $H(t)\rightarrow \infty 
$ when $t\rightarrow t_{1}$ for some 
\[
t_{1}\leq H(0)\left[ \nu H^{\left( 1\right) }(0)\right] ^{-1}. 
\]%
We rewrite the energy identity as%
\[
E\left( t\right) =\left\Vert Bu_{t}\right\Vert ^{2}+\left( \left[
B^{2}\left( A\ast u-\left( a\ast \Delta \right) \right) u\right] ,u\right)
+\left( f\left( u\right) ,u\right) =E\left( 0\right) . 
\]%
We prove here the following

\textbf{Theorem 5.1. }Assume the Condition 4.1 is satisfied and $s+r\geq 1$.
Let $B\varphi $, $B\psi \in L^{2}$. If there are some positive numbers $\nu $%
, $t_{0}$ and $b$ such that

\bigskip 
\[
\left( 1+2\nu \right) b+d\leq -E\left( 0\right) \text{, }4b\left( 1+\nu
\right) \left( t+t_{0}\right) \leq 2b-2E\left( 0\right) 
\]%
and 
\[
4\nu \left\Vert Bu\right\Vert ^{2}\left\Vert Bu_{t}\right\Vert ^{2}\leq
\varkappa _{1}\left\Vert Bu\right\Vert ^{2}+\varkappa _{2}\left\Vert
Bu_{t}\right\Vert ^{2}+\phi \left( t\right) , 
\]%
for the solution $u\in C^{\left( 2,s\right) }\left( A,E\right) $ of $%
(1.1)-(1.2)$ and%
\begin{equation}
E\left( 0\right) =\left\Vert B\psi \right\Vert ^{2}+\left( B^{2}\left[ A\ast
u-a\ast \Delta u\right] ,u\right) +\left( f\left( \varphi \right) ,\varphi
\right) <0\text{,}  \tag{5.1}
\end{equation}%
for all $t\geq 0$, where 
\[
\left( 1+2\nu \right) b\leq -E\left( 0\right) \text{,} 
\]%
\[
\varkappa _{1}=\left[ 2b-2E\left( 0\right) -4b\left( 1+\nu \right) \left(
t+t_{0}\right) \right] , 
\]%
\[
\varkappa _{2}=4b\left( t+t_{0}\right) \left[ \left( t+t_{0}\right) -\left(
1+\nu \right) \right] \text{, } 
\]%
\[
\phi \left( t\right) =2b\left[ -E\left( 0\right) -\left( 1+2\nu \right) b%
\right] \left( t+t_{0}\right) ^{2}. 
\]%
Then the solution $u$ blows up in finite time.

\textbf{Proof. }Assume that there is a global solution. Then $Bu(t)$, $%
Bu_{t}(t)\in $ $L^{2}$ for all $t>0$. Let 
\[
H\left( t\right) =\left\Vert Bu\right\Vert ^{2}+b\left( t+t_{0}\right) ^{2}%
\text{.} 
\]%
for some positive $b$ and $t_{0}$ that will be determined later. We have \ 
\begin{equation}
H^{\left( 1\right) }\left( t\right) =2\left( Bu,Bu_{t}\right) +2b\left(
t+t_{0}\right) ,  \tag{5.2}
\end{equation}

\[
H^{\left( 2\right) }\left( t\right) =2\left\Vert Bu_{t}\right\Vert
^{2}+2\left( Bu,Bu_{tt}\right) +2b. 
\]

From $\left( 1.1\right) $ and $\left( 5.1\right) $ we get \ 

\[
\left( Bu,Bu_{tt}\right) =\left( u,B^{2}u_{tt}\right) = 
\]%
\[
\left( u,B^{2}\left[ a\ast \Delta u-A\ast u+\Delta g\ast f\left( u\right) %
\right] \right) = 
\]%
\begin{equation}
\left[ \left( u,B^{2}\left( a\ast \Delta u\right) \right) -\left(
u,B^{2}A\ast u\right) -\left( u,f\left( u\right) \right) \right] =  \tag{5.3}
\end{equation}%
\[
\left\Vert Bu_{t}\right\Vert ^{2}-E\left( 0\right) . 
\]

From $\left( 5.2\right) $ and $\left( 5.3\right) $, we obtain%
\begin{equation}
H^{\left( 2\right) }\left( t\right) \geq 4\left\Vert Bu_{t}\right\Vert
^{2}-2E\left( 0\right) +2b.  \tag{5.4}
\end{equation}

\bigskip On the other hand, in view of Cauchy-Schwartz inequality, we have%
\[
\left( H^{\left( 1\right) }\left( t\right) \right) ^{2}=\left[ 2\left(
Bu,Bu_{t}\right) +2b\left( t+t_{0}\right) \right] ^{2}\leq 
\]%
\[
4\left[ \left\Vert Bu\right\Vert ^{2}\left\Vert Bu_{t}\right\Vert
^{2}+b\left( t+t_{0}\right) \left( \left\Vert Bu\right\Vert ^{2}+\left\Vert
Bu_{t}\right\Vert ^{2}\right) \right] + 
\]%
\begin{equation}
4b^{2}\left( t+t_{0}\right) ^{2}.  \tag{5.5}
\end{equation}

Hence, by $\left( 5.2\right) $, $\left( 5.4\right) $ and $\left( 5.5\right) $%
,\ we obtain 
\[
H^{\left( 2\right) }H-\left( 1+\nu \right) \left( H^{\left( 1\right)
}\right) ^{2}\geq 
\]%
\[
\left[ 4\left\Vert Bu_{t}\right\Vert ^{2}+2b-2E\left( 0\right) \right] \left[
\left\Vert Bu\right\Vert ^{2}+b\left( t+t_{0}\right) ^{2}\right] - 
\]%
\[
4\left( 1+\nu \right) \left[ \left\Vert Bu\right\Vert ^{2}\left\Vert
Bu_{t}\right\Vert ^{2}+b\left( t+t_{0}\right) \left( \left\Vert
Bu\right\Vert ^{2}+\left\Vert Bu_{t}\right\Vert ^{2}\right) \right] - 
\]%
\[
4\left( 1+\nu \right) b^{2}\left( t+t_{0}\right) ^{2}\geq 
\]%
\[
2b\left[ -E\left( 0\right) -\left( 1+2\nu \right) b\right] \left(
t+t_{0}\right) ^{2}+ 
\]

\[
\left[ 2b-2E\left( 0\right) -4b\left( 1+\nu \right) \left( t+t_{0}\right) %
\right] \left\Vert Bu\right\Vert ^{2}+ 
\]%
\[
4b\left( t+t_{0}\right) \left[ \left( t+t_{0}\right) -\left( 1+\nu \right) %
\right] \left\Vert Bu_{t}\right\Vert ^{2}-4\nu \left\Vert Bu\right\Vert
^{2}\left\Vert Bu_{t}\right\Vert ^{2}\geq 0, 
\]%
when%
\[
\left( 1+2\nu \right) b\leq -E\left( 0\right) \text{, }4b\left( 1+\nu
\right) \left( t+t_{0}\right) \leq 2b-2E\left( 0\right) 
\]%
and 
\[
4\nu \left\Vert Bu\right\Vert ^{2}\left\Vert Bu_{t}\right\Vert ^{2}\leq
\varkappa _{1}\left\Vert Bu\right\Vert ^{2}+\varkappa _{2}\left\Vert
Bu_{t}\right\Vert ^{2}+\phi \left( t\right) , 
\]%
for all $t\geq 0$, where%
\[
\left( 1+2\nu \right) b\leq -E\left( 0\right) \text{,} 
\]%
\[
\varkappa _{1}=\left[ 2b-2E\left( 0\right) -4b\left( 1+\nu \right) \left(
t+t_{0}\right) \right] , 
\]%
\[
\varkappa _{2}=4b\left( t+t_{0}\right) \left[ \left( t+t_{0}\right) -\left(
1+\nu \right) \right] \text{, } 
\]%
\[
\phi \left( t\right) =2b\left[ -E\left( 0\right) -\left( 1+2\nu \right) b%
\right] \left( t+t_{0}\right) ^{2}. 
\]

Then by Theorem 5.1 we obtain the assertion.

\begin{center}
\textbf{\ \ 6. Applications}

\textbf{6.1}.\textbf{The Cauchy problem for the system of nonlocal WEs }
\end{center}

\bigskip Consider the problem $\left( 1.3\right) $. Let%
\[
\text{ }l_{q}=\left\{ \text{ }u=\left\{ u_{j}\right\} ,\text{ }j=1,2,...N%
\text{, }\left\Vert u\right\Vert _{l_{q}}=\left( \sum\limits_{j=1}^{\infty
}\left\vert u_{j}\right\vert ^{q}\right) ^{\frac{1}{q}}<\infty \right\} ,
\]%
(see $\left[ \text{23, \S\ 1.18}\right] $. Let $A_{1}$ be the operator in $%
l_{p}$ defined by%
\[
\text{ }A_{1}=\left[ a_{jm}\left( x\right) \right] ,\text{ }%
a_{jm}=b_{j}\left( x\right) 2^{\sigma m}\text{, }m,j=1,2,...\infty \text{, }%
D\left( A_{1}\right) =\text{ }l_{q}^{\sigma }=
\]

\[
\left\{ \text{ }u=\left\{ u_{j}\right\} ,\text{ }j=1,2,...\infty \text{, }%
\left\Vert u\right\Vert _{l_{q}^{\sigma }}=\left( \sum\limits_{j=1}^{\infty
}2^{\sigma j}\left\vert u_{j}\right\vert ^{q}\right) ^{\frac{1}{q}}<\infty
\right\} \text{, }\sigma >0.
\]

Let \ 
\[
W^{s,p}\left( E\right) =W^{s,p}\left( \mathbb{R}^{n};E\right) ,\text{ }%
W^{s}\left( E\right) =W^{s,2}\left( \mathbb{R}^{n};E\right) , 
\]%
\[
Y^{s,p,\sigma }=W^{s,p}\left( \mathbb{R}^{n};l_{q}\right) \cap L^{p}\left( 
\mathbb{R}^{n};l_{q}^{\sigma }\right) ,\text{ }1\leq q\leq \infty , 
\]%
\[
W_{0}\left( l_{q}\right) =W^{s\left( 1-\frac{1}{2p}\right) ,p}\left( \mathbb{%
R}^{n};l_{q}\right) \cap L^{p}\left( \mathbb{R}^{n};l_{q}^{\sigma \left( 1-%
\frac{1}{2p}\right) }\right) . 
\]

\ Let $f=\left\{ f_{m}\right\} $, $m=1,2,...\infty $ and 
\[
\eta _{1}=\eta _{1}\left( \xi \right) =\left[ \hat{a}\left( \xi \right)
\left\vert \xi \right\vert ^{2}+\hat{A}_{1}\left( \xi \right) \right] ^{%
\frac{1}{2}}.
\]

Here, 
\[
E_{ip}\left( l_{q}\right) =W^{s\left( 1-\theta _{i}\right) ,p}\left( \mathbb{%
R}^{n};l_{q}\right) \cap L^{p}\left( \mathbb{R}^{n};l_{q}^{\sigma \left(
1-\theta _{i}\right) }\right) , 
\]%
where 
\[
\theta _{j}=\frac{1+ip}{2p}\text{, }i=0,1\text{.} 
\]

\ From Theorem 3.1 we obtain the following result \ 

\textbf{Theorem 6.1.} Assume: (1) $0<\alpha <1-\frac{1}{2p}$, $\varphi \in
E_{0p}\left( l_{q}\right) $, $\psi \in $ $E_{1p}\left( l_{q}\right) $ and $%
s>1+\frac{n}{p}$ for $p\in \left[ 1,\infty \right] $, $q\in \left( 1,\infty
\right) $; (2) the assumptions (1)-(2) of Condition 2.1 are satisfied; (3) $%
\hat{b}_{j}=b_{j}\left( \xi \right) $ are nonnegat\i ve bounded d\i
fferentiable functions on $\mathbb{R}^{n}$ and $a+\hat{b}_{j}\left( \xi
\right) >0$ for $\xi \in \mathbb{R}^{n}$, $D^{\alpha }\hat{b}_{j}$ are
uniformly bounded on $\mathbb{R}^{n}$ for $\alpha =\left( \alpha _{1},\alpha
_{2},...,\alpha _{n}\right) $, $\left\vert \alpha \right\vert \leq n$ and
the uniform estimate holds 
\[
\dsum\limits_{j=1}^{N}\left\vert D^{\alpha }\hat{b}_{j}\left( \xi \right)
\right\vert ^{2}\left[ \hat{a}\left( \xi \right) \left\vert \xi \right\vert
^{2}+\hat{b}_{j}\left( \xi \right) \right] ^{-1}\leq M;
\]%
(4) the kernel $g_{mj}$ are bounded integrable functions, whose Fourier
transform satisfies%
\[
0\leq \dsum\limits_{j=m,j}^{N}\left\vert \hat{g}_{mj}\left( \xi \right)
\right\vert ^{2}\lesssim \left( 1+\left\vert \xi \right\vert ^{2}\right) ^{-%
\frac{r}{2}}\text{ for all }\xi \in \mathbb{R}^{n}\text{ and }r\geq 2;
\]%
(5) the function%
\[
u\rightarrow f\left( x,t,u\right) :\mathbb{R}^{n}\times \left[ 0,T\right]
\times W_{0}\left( l_{q}\right) \rightarrow l_{q}
\]%
is a measurable in $\left( x,t\right) \in \mathbb{R}^{n}\times \left[ 0,T%
\right] $ for $u\in W_{0}\left( l_{q}\right) $; Moreover, $f\left(
x,t,u\right) $ is continuous in $u\in W_{0}\left( l_{q}\right) $ and $f\in
C^{\left[ s\right] +1}\left( W_{0}\left( l_{q}\right) ;l_{q}\right) $
uniformly in $x\in \mathbb{R}^{n},$ $t\in \left[ 0,T\right] $. Then problem $%
\left( 1.3\right) $ has a unique local strange solution%
\[
u\in C^{\left( 2\right) }\left( \left[ 0,\right. \left. T_{0}\right)
;Y_{\infty }^{s,p}\left( A_{1},l_{q}\right) \right) ,
\]%
where $T_{0}$ is a maximal time interval that is appropriately small
relative to $M$. Moreover, if

\[
\sup_{t\in \left[ 0\right. ,\left. T_{0}\right) }\left( \left\Vert
u\right\Vert _{Y_{\infty }^{s,p}\left( \hat{A}_{1}^{\alpha };l_{q}\right)
}+\left\Vert u_{t}\right\Vert _{Y_{\infty }^{s,p}\left( \hat{A}_{1}^{\alpha
};l_{q}\right) }\right) <\infty 
\]%
then $T_{0}=\infty .$

\ \textbf{Proof. }It is known that $L^{p}\left( \mathbb{R}^{n};l_{q}\right) $
is a UMD space for $p,q\in \left( 1,\infty \right) $ (see e.g $\left[ 25%
\right] $). By Remark 2.1, by definition of $W^{s,p}\left(
A_{1},l_{q}\right) $ and by real interpolation of Banach spaces (see e.g. $%
\left[ \text{23, \S 1.3, 1.18}\right] $), we have 
\[
\text{ }\mathbb{E}_{ip}=\left( W^{s,,p}\left( \mathbb{R}^{n};l_{q}^{\sigma
},l_{q}\right) ,L_{p}\left( \mathbb{R}^{n};l_{q}\right) _{\theta
_{i},p}\right) =W^{s\left( 1-\theta _{i}\right) ,p}\left( \mathbb{R}%
^{n};l_{q}^{\sigma \left( 1-\theta _{i}\right) },l_{q}\right) = 
\]%
\[
W^{s\left( 1-\theta _{i}\right) ,p}\left( \mathbb{R}^{n};l_{q}\right) \cap
L^{p}\left( \mathbb{R}^{n};l_{q}^{\sigma \left( 1-\theta _{i}\right)
}\right) =E_{0i}\left( l_{q}\right) \text{, }i=0,1. 
\]%
By assumptions (1), (2) we obtain that $\hat{A}_{1}\left( \xi \right) $ is
uniformly sectorial in $l_{q}$, $\hat{A}_{1}\left( \xi \right) \in \sigma
\left( M_{0},\omega ,l_{q}\right) ,$ $\eta _{1}\left( \xi \right) \neq 0$\
for all $\ \xi \in \mathbb{R}^{n}$\ and 
\[
\left\Vert D^{\alpha }\hat{A}_{1}\left( \xi \right) \eta _{1}^{-1}\left( \xi
\right) \right\Vert _{B\left( l_{q}\right) }\leq M 
\]%
for $\alpha =\left( \alpha _{1},\alpha _{2},...,\alpha _{n}\right) $, $%
\left\vert \alpha \right\vert \leq n$. Hence, by (4), (5), all conditions of
Theorem 3.2 are hold, i,e., we get the conclusion.\ 

Let $G$ be a function defined by $\left( 4.15\right) $.

\textbf{Theorem 6.2.} Assume: (a) (1)-(3) assumptions of Theorem 6.1 are
satisfied for $p=2$ and 
\[
\left\Vert \hat{g}\left( \xi \right) \right\Vert _{B\left( l_{2}\right)
}\lesssim \left( 1+\left\vert \xi \right\vert ^{2}\right) ^{-\frac{r}{2}}%
\text{ for }r\leq 2\left( s+1\right) , 
\]%
\[
\left\Vert \hat{g}^{\frac{1}{2}}\left( \xi \right) \right\Vert _{B\left(
l_{2}\right) }\lesssim \left\vert \xi \right\vert \left( 1+\left\vert \xi
\right\vert ^{2}\right) ^{\frac{s}{2}}\text{ for all }\xi \in \mathbb{R}%
^{n}; 
\]%
(c) $f_{m}\in C^{\left[ s\right] }\left( \mathbb{R};l_{2}\right) $ with $%
f(0)=0$ and%
\[
\dsum\limits_{m=1}^{N}\left\vert \hat{f}_{m}\left( u\right) \left( \xi
\right) \right\vert ^{2}<\infty \text{ for all }u=\left\{ u_{m}\right\} \in
C^{\left( 2\right) }\left( \left[ 0,\right. \left. \infty \right) ;Y_{\infty
}^{s,2}\left( A_{1};l_{2}\right) \right) ; 
\]

(d) $B\varphi $, $B\psi \in L^{2}\left( \mathbb{R}^{n};l_{2}\right) $ and $%
\left( \varphi \text{, }f\left( \varphi \right) \right) \in L^{2}\left( 
\mathbb{R}^{n};l_{2}\right) $; (e) there is some $k>0$ so that  
\[
\left( \varphi \text{, }f\left( \varphi \right) \right) _{L^{2}\left( 
\mathbb{R}^{n};l_{2}\right) }\geq -k\left\Vert \text{ }\varphi \right\Vert
_{L^{2}\left( \mathbb{R}^{n};l_{2}\right) }.
\]
Then there is some $T>0$ such that problem $(1.3)$ has a global solution%
\[
u\in C^{\left( 2\right) }\left( \left[ 0,\right. \left. \infty \right)
;Y_{\infty }^{s,2}\left( A_{1};l_{2}\right) \right) .
\]

\ \textbf{Proof. }From the assumptions (a), (b) it is clear to see that the
Condition 4.1 holds for $H=l_{2}$ and $r>2+\frac{n}{2}$. By (c), (d), (e)
all other assumptions of Theorem 4.1 are satisfied. Hence, we obtain the
assertion.

\begin{center}
\textbf{6.2.} \textbf{The mixed problem for degenerate nonlocal WE}
\end{center}

Consider the problem $\left( 1.5\right) -\left( 1.7\right) $. Let%
\[
Y^{s,p,2}=W^{s,p}\left( \mathbb{R}^{n};L^{p_{1}}\left( 0,1\right) \right)
\cap L^{p}\left( \mathbb{R}^{n};W^{\left[ 2\right] ,p_{1}}\left( 0,1\right)
\right) \text{, }1\leq p\leq \infty , 
\]

Let $A_{2}$ is the operator in $L^{p_{1}}\left( 0,1\right) $ defined by $%
\left( 1.4\right) $ and let%
\[
\eta _{2}=\eta _{2}\left( \xi \right) =\left[ a\left\vert \xi \right\vert
^{2}+\hat{A}_{2}\left( \xi \right) \right] ^{\frac{1}{2}}. 
\]

Here,%
\[
E_{ip}\left( L^{p_{1}}\right) =W^{\left[ s\left( 1-\theta _{i}\right) \right]
,p}\left( \mathbb{R}^{n};L^{p_{1}}\left( 0,1\right) \right) \cap L^{p}\left( 
\mathbb{R}^{n};W^{\left[ 2\left( 1-\theta _{i}\right) \right] ,p_{1}}\left(
0,1\right) \right) \text{,} 
\]%
where 
\[
\theta _{i}=\frac{1+ip}{2p}\text{, }i=0,1. 
\]

Now, we present the following result:

\textbf{Condition 6.1 }Assume;

(1) $0\leq \gamma <\frac{1}{p_{1}}$ for $p_{1}\in \left( 1,\infty \right) $
and $\alpha _{1}\beta _{2}-\alpha _{2}\beta _{1}\neq 0;$

(2) $0<\alpha <1-\frac{1}{2p}$, $\varphi \in E_{0p}\left( L^{p_{1}}\right) $%
, $\psi \in $ $E_{1p}\left( L^{p_{1}}\right) $ and $s>1+\frac{n}{p}$ for $%
p\in \left[ 1,\infty \right] $, $p_{1}\in \left( 1,\infty \right) $;

(2) $b_{1}$ and $b_{2}$ are complex valued functions on $\left( 0,1\right) $%
. Morover, $b_{1}\in C\left[ 0,1\right] ,$ $b_{1}\left( 0\right)
=b_{1}\left( 1\right) $, $b_{2}\in L_{\infty }\left( 0,1\right) $ and $%
\left\vert b_{2}\left( x\right) \right\vert \leq C$ $\left\vert b_{1}^{\frac{%
1}{2}-\mu }\left( x\right) \right\vert $ for $0<\mu <\frac{1}{2}$ and for
a.a. $x\in \left( 0,1\right) ;$

(3) the assumptions (1)-(2) of Condition 2.1 are satisfied; $D^{\alpha }\hat{%
b}_{j}$, $j=1,$ $2$ are uniformly bounded on $\mathbb{R}^{n}$ for all $%
\alpha =\left( \alpha _{1},\alpha _{2},...,\alpha _{n}\right) $ with $%
\left\vert \alpha \right\vert \leq n$ and $\eta _{2}\left( \xi \right) \neq 0
$\ for all $\xi \in \mathbb{R}^{n};$

(4)$\ $for $\alpha =\left( \alpha _{1},\alpha _{2},...,\alpha _{n}\right) ,$ 
$\left\vert \alpha \right\vert \leq n$ the uniform estimate holds 
\[
\left\Vert \left[ D^{\alpha }\hat{A}_{2}\left( \xi \right) \right] \eta
_{2}^{-1}\left( \xi \right) \right\Vert _{B\left( L^{p_{1}}\left( 0,1\right)
\right) }\leq M. 
\]

(5) the function%
\[
u\rightarrow f\left( x,t,u\right) :\mathbb{R}^{n}\times \left[ 0,T\right]
\times W_{0}\left( L^{p_{1}}\left( 0,1\right) \right) \rightarrow
L^{p_{1}}\left( 0,1\right) 
\]%
is a measurable in $\left( x,t\right) \in \mathbb{R}^{n}\times \left[ 0,T%
\right] $ for $u\in W_{0}\left( L^{p_{1}}\left( 0,1\right) \right) $; $%
f\left( x,t,u\right) $. Moreover, $f\left( x,t,u\right) $ is continuous in $%
u\in W_{0}\left( L^{p_{1}}\left( 0,1\right) \right) $ and%
\[
f\left( x,t,u\right) \in C^{\left[ s\right] +1}\left( W_{0}\left(
L^{p_{1}}\left( 0,1\right) \right) ;L^{p_{1}}\left( 0,1\right) \right) 
\]%
uniformly with respect to $x\in \mathbb{R}^{n}$, $t\in \left[ 0,T\right] .$

\bigskip \textbf{Theorem 6.3.} Assume that the Condition 6.1 is satisfied.
Then problem $\left( 1.5\right) -\left( 1.7\right) $ has a unique local
strange solution 
\[
u\in C^{\left( 2\right) }\left( \left[ 0,\right. \left. T_{0}\right)
;Y_{\infty }^{s,p}\left( A_{2},L^{p_{1}}\left( 0,1\right) \right) \right) , 
\]%
where $T_{0}$ is a maximal time interval that is appropriately small
relative to $M$. Moreover, if

\[
\sup_{t\in \left[ 0\right. ,\left. T_{0}\right) }\left( \left\Vert
u\right\Vert _{Y_{\infty }^{s,p}\left( \hat{A}_{2}^{\alpha };L^{2}\left(
0,1\right) \right) }+\left\Vert u_{t}\right\Vert _{Y_{\infty }^{s,p}\left( 
\hat{A}_{2}^{\alpha };L^{2}\left( 0,1\right) \right) }\right) <\infty 
\]%
then $T_{0}=\infty .$

\ \textbf{Proof.}\ It is known (see e.g. $\left[ 13\right] $) that $%
L^{p_{1}}\left( 0,1\right) $ is a UMD space for $p_{1}\in \left( 1,\infty
\right) $. By definition of $W^{s,p}\left( A_{2},L^{p_{1}}\left( 0,1\right)
\right) $ and by real interpolation of Banach spaces (see e.g. $\left[ \text{%
23, \S 1.3}\right] $) we have

\[
\text{ }\mathbb{E}_{ip}=W^{s,,p}\left( \mathbb{R}^{n};W^{\left[ 2\right]
,,p_{1}}\left( 0,1\right) ,L^{p_{1}}\left( 0,1\right) ,L^{p}\mathbb{R}%
^{n};L^{p_{1}}\left( 0,1\right) \right) _{\theta _{i},p}=
\]%
\[
W^{s\left( 1-\theta _{i}\right) ,p}\left( \mathbb{R}^{n};W^{\left[ 2\left(
1-\theta _{i}\right) \right] ,p_{1}}\left( 0,1\right) ,L^{p_{1}}\left(
0,1\right) \right) =E_{ip}\left( L^{p_{1}}\right) .
\]%
In view of $\left[ \text{26, Theorem 4.1}\right] $ we obtain that $\hat{A}%
_{2}\left( \xi \right) $ is uniformly sectorial in $L^{p_{1}}\left(
0,1\right) $ and 
\[
\hat{A}_{2}\left( \xi \right) \in \sigma \left( M_{0},\omega
,L^{p_{1}}\left( 0,1\right) \right) .
\]%
Moreover, by using the assumptions (1), (2) we deduced that $\eta _{2}\left(
\xi \right) \neq 0$\ for all $\ \xi \in \mathbb{R}^{n}$\ and 
\[
\left\Vert D^{\alpha }\hat{A}_{2}\left( \xi \right) \eta _{2}^{-1}\left( \xi
\right) \right\Vert _{B\left( L^{p_{1}}\left( 0,1\right) \right) }\leq M.
\]%
for $\alpha =\left( \alpha _{1},\alpha _{2},...,\alpha _{n}\right) ,$ $%
\left\vert \alpha \right\vert \leq n$. Hence, by hypothesis (3), (4) of the
Condition 5.1 we get that all hypothesis of Theorem 3.2 are hold, i,e., we
obtain the conclusion.\ Let $G=\left( 0,1\right) \times \mathbb{R}^{n}$

\bigskip \textbf{Theorem 6.4.} Assume the Condition 6.1 is satisfied for $%
p_{1}=2$.\ Suppose $f\in C^{\left[ s\right] }\left( \mathbb{R};L^{2}\left(
\left( 0,T\right) \right) \right) $ with $f(0)=0$. Let the kernel $g_{mj}$
be bounded integrable functions and

\[
\left\Vert \hat{g}\left( \xi \right) \right\Vert _{B\left( l_{2}\right)
}\lesssim \left( 1+\left\vert \xi \right\vert ^{2}\right) ^{-\frac{r}{2}}%
\text{ for }r\leq 2\left( s+1\right) , 
\]%
\[
\left\Vert \hat{g}^{\frac{1}{2}}\left( \xi \right) \right\Vert _{B\left(
l_{2}\right) }\lesssim \left\vert \xi \right\vert \left( 1+\left\vert \xi
\right\vert ^{2}\right) ^{\frac{s}{2}}\text{ for all }\xi \in \mathbb{R}^{n} 
\]

Moreover, let $\ B\varphi ,$ $B\psi \in L^{2}\left( G\right) $, and $\left(
\varphi \text{, }f\left( \varphi \right) \right) \in L^{2}\left( G\right) $;
(e) there is some $k>0$ so that 
\[
\left( \varphi \text{, }f\left( \varphi \right) \right) _{L^{2}\left( 
\mathbb{G}\right) }\geq -k\left\Vert \text{ }\varphi \right\Vert
_{L^{2}\left( G\right) }.
\]

Then there is some $T>0$ such that the problem $(1.5)-\left( 1.7\right) $
has a global solution%
\[
u\in C^{2}\left( \left[ 0,\right. \left. \infty \right) ;Y_{\infty
}^{s,2}\right) .
\]

\textbf{Proof.} Indeed, by assumptions all conditions of Theorem 4.1. are
satisfied for $H=L^{2}\left( 0,1\right) $, i.e. we obtain the assertion.

\begin{enumerate}
\item \textbf{References}
\end{enumerate}

\begin{quote}
\ \ \ \ \ \ \ \ \ \ \ \ \ \ \ \ \ \ \ \ \ \ \ \ 
\end{quote}

\begin{enumerate}
\item M. Arndt and M. Griebel, Derivation of higher order gradient continuum
models from atomistic models for crystalline solids Multiscale Modeling
Simul. (2005)4, 531--62.

\item A. Ashyralyev, N. Aggez, Nonlocal boundary value hyperbolic problems
involving Integral conditions, Bound.Value Probl., 2014 V (2014):214.

\item X. Blanc, C. LeBris, P. L. Lions, Atomistic to continuum limits for
computational materials science, ESAIM- Math. Modelling Numer. Anal.
(2007)41, 391--426.

\item J.L. Bona, R.L. Sachs, Global existence of smooth solutions and
stability of solitary waves for a generalized Boussinesq equation, Comm.
Math. Phys. 118 (1988), 15--29.

\item A. Constantin and L. Molinet, The initial value problem for a
generalized Boussinesq equation, Diff.Integral Eqns. (2002)15, 1061--72.

\item G. Chen and S. Wang, Existence and nonexistence of global solutions
for the generalized IMBq equation Nonlinear Anal.---Theory Methods Appl.
(1999)36, 961--80.

\item R. Coifman and Y. Meyer, Wavelets. Calder%
%TCIMACRO{\U{b4}}%
%BeginExpansion
\'{}%
%EndExpansion
on-Zygmund and multilinear operators, Cambridge University Press, 1997.

\item M. Dafermos, I. Rodnianski, Y. Shlapentokh-Rothman, Decay for
solutions of the wave equation on Kerr exterior spacetimes III: The full
subextremalcase \TEXTsymbol{\vert}a\TEXTsymbol{\vert}\TEXTsymbol{<}M, Anal.
Math, 183 (2016),787-913.

\item A. De Godefroy, Blow up of solutions of a generalized Boussinesq
equation IMA J. Appl. Math.(1998) 60 123--38.

\item N. Duruk, H.A. Erbay and A. Erkip, Global existence and blow-up for a
class of nonlocal nonlinear Cauchy problems arising in elasticity,
Nonlinearity, (2010)23, 107--118.

\item A. C. Eringen, Nonlocal Continuum Field Theories, New York, Springer
(2002).

\item H. O. Fattorini, Second order linear differential equations in Banach
spaces, in North Holland Mathematics Studies,\ V. 108, North-Holland,
Amsterdam, 1985.

\item M. Girardi, L. Weis, Operator-valued Fourier multiplier theorems on $%
L_{p}$($X$) and geometry of Banach spaces, J. Funct. Anal., 204(2) (2003),
320--354.

\item Z. Huang, Formulations of nonlocal continuum mechanics based on a new
definition of stress tensor Acta Mech. (2006)187, 11--27.

\item T. Kato, G. Ponce, Commutator estimates and the Euler and
Navier--Stokes equations, Comm. Pure Appl. Math. (1988)41, 891--907.

\item V. K. Kalantarov and O. A. Ladyzhenskaya, The occurence of collapse
for quasilinear equation of parabolic and hyperbolic types Journal of Soviet
Mathematics (10) (1978) 53-70.

\item M. Lazar, G. A. Maugin and E. C. Aifantis, On a theory of nonlocal
elasticity of bi-Helmholtz type and some applications Int. J. Solids and
Struct. (2006)43, 1404--21.

\item F. Linares, Global existence of small solutions for a generalized
Boussinesq equation, J. Differential Equations 106 (1993), 257--293.

\item Y. Liu, Instability and blow-up of solutions to a generalized
Boussinesq equation, SIAM J. Math. Anal. 26 (1995), 1527--1546.

\item V.G. Makhankov, Dynamics of classical solutions (in non-integrable
systems), Phys. Lett. C 35(1978), 1--128.

\item C. Polizzotto, Nonlocal elasticity and related variational principles
Int. J. Solids Struct. ( 2001) 38 7359--80.

\item A. Pazy, Semigroups of linear operators and applications to partial
differential equations. Springer, Berlin, 1983.

\item L. S. Pulkina, A non local problem with integral conditions for
hyperbolice quations, Electron. J. Differ. Equ.(1999)45, 1-6.

\item C. A. Silling, Reformulation of elasticity theory for discontinuities
and long-range forces J. Mech. Phys. Solids (2000)48 175-209.

\item V. B. Shakhmurov, Embedding and separable differential operators in
Sobolev-Lions type spaces, Math. Notes, 84(2008) (6), 906-926.

\item V. B. Shakhmurov, Linear and nonlinear abstract differential equations
with small parameters, Banach J. Math. Anal. 10 (2016)(1), 147--168.

\item H. Triebel, Interpolation theory, Function spaces, Differential
operators, North-Holland, Amsterdam, 1978.

\item H. Triebel, Fractals and spectra, Birkhauser Verlag, Related to
Fourier analysis and function spaces, Basel, 1997.

\item G.B. Whitham, Linear and Nonlinear Waves, Wiley--Interscience, New
York, 1975.

\item S. Wang, G. Chen, Small amplitude solutions of the generalized IMBq
equation, J. Math. Anal. Appl. 274 (2002) 846--866.

\item S.Wang and G.Chen, Cauchy problem of the generalized double dispersion
equation Nonlinear Anal. Theory Methods Appl. (2006 )64 159--73.

\item Ta-Tsien Li, Y. Jinsemi, Global $C^{1}$ solution to the mixed initial
boundary value problems for quasilinear hyperbolic system, Chinese Ann.
Math. (22)03 (2001), 325-336.

\item N.J. Zabusky, Nonlinear Partial Differential Equations, Academic
Press, New York, 1967.
\end{enumerate}

\end{document}